\documentclass[12pt]{article}
\usepackage[backref=page,colorlinks,bookmarksopen=true]{hyperref}
\usepackage{amssymb}
\usepackage{latexsym}
\usepackage{amsmath,mathabx}
\usepackage[all]{xy}
\newtheorem{theorem}{Theorem}[section]
\newtheorem{lemma}[theorem]{Lemma}
\newtheorem{proposition}[theorem]{Proposition}
\newtheorem{corollary}[theorem]{Corollary}
\newtheorem{definition}[theorem]{Definition}
\newtheorem{example}[theorem]{Example}

\newtheorem{remark}[theorem]{Remark}

\begin{document}
\title{Tannaka-Krein reconstruction for coactions of finite quantum groupoids}
\author{Leonid Vainerman Jean-Michel Vallin}
\date{16 June 2016}
\maketitle

\begin{abstract} We study coactions of finite quantum groupoids on unital $C^*$-algebras and obtain the
Tannaka-Krein reconstruction theorem for them.
\end{abstract}

\tableofcontents

\footnote {AMS Subject Classification [2010]{:  18D10, 16T05., 46L05 }}

\footnote{Keywords {: Coactions and corepresentations of quantum groupoids, $C^*$-categories, reconstruction theorem.}}
\newpage
\newenvironment{dm}{\hspace*{0,15in} {{\bf Proof.}}}{$\square$}
%%%%%%%%%%%%%%%%%%%%%%%%%%%%%%%%%%%%%%%%%%%%%%%%%%%%%%%%%%%%%%%%%%%%%%%%%%%%%%%%%%%%%%%%%%%%%%%%%%%%%%%%%%%%%%
\begin{section}{Introduction}
As shown in \cite{NV1},
\cite{NV2}, weak Hopf $C^*$-algebra in the sense of \cite{BNSz} and their coideal $C^*$-subalgebras play important role in the
description of Jones's tower of $II_1$-subfactors with finite index and finite depth.  It is also known that any fusion category
can be realized as a representation category of a weak Hopf algebra \cite{EGNO}. This explains the interest in the construction
of concrete examples of these objects and in the classification of their coideal subalgebras.

In this paper we use the term "a finite quantum groupoid" instead of "a weak Hopf $C^*$-algebra" because a groupoid $C^*$-algebra
and an algebra of functions on a usual finite groupoid carry this type of a structure. Some particular constructions of finite
quantum groupoids were proposed in \cite {Val2} and \cite{NV5} (see also the survey \cite{NV} and references therein), but the
general way to construct them is the application of the Hayashi's reconstruction theorem of Tannaka-Krein type \cite{Ha} to concrete
tensor categories.

The above approach was used in \cite{M}, where a series of concrete finite quantum groupoids was constructed using  Tambara-Yamagami
categories \cite{TY}. These categories belong to the much wider family of ${\mathbb Z}/2{\mathbb Z}$-extensions of pointed fusion
categories classified in \cite{VV}. The problem of the description of coideal $C^*$-subalgebras of a given finite quantum groupoid is
even harder, and until now only two concrete families of such subalgebras constructed in \cite{M} "by hand" are known.

On the other hand, a similar problem exists in the theory of compact quantum groups (CQG), where a lot of concrete examples were
constructed using well-known Woronowicz's reconstruction theorem of Tannaka-Krein type \cite{Wor}. As for the description of coideal
$C^*$-subalgebras of a given CQG, the recent papers on categorical duality for (co)actions of CQG on $C^*$-algebras proved to be
very useful - see \cite{CY}, \cite{CY1}, \cite{Nes1}, \cite{NY1}.

Our aim, here and in a subsequent work, is to develop a similar approach for the study of finite quantum groupoid coactions
on $C^*$-algebras and to apply it to the description of coideal $C^*$-subalgebras of concrete finite quantum groupoids. The
first step in this project is to formulate and to prove the general duality statement - Theorem 1.1 below.

Let us describe the structure of the paper. In Section 2 we recall basic definitions and results on finite quantum groupoids
following \cite{BNSz} and \cite{NV}. We also translate the representation theory of these objects treated in \cite{BoSz} and
\cite{NTV} into the language of unitary corepresentations and $C^*$-tensor categories suitable for the construction of the
categorical duality. Finally, we translate into this language the reconstruction theorem proved in \cite{Ha} and \cite{SZ}.

In Section 3 we develop the theory parallel to the one of CQG coactions \cite{Boca}. Doing this, we simplify significantly, in our
particular case, some constructions related to coactions of general measured quantum groupoids - see \cite{E2}, \cite{Val2}. Let
$\mathfrak a$ be a coaction of a finite quantum groupoid $\mathfrak G$ on a unital $C^*$-algebra $A$ (called a
$\mathfrak G$-$C^*$-algebra) we get the canonical implementation of $\mathfrak a$ and study the properties of the spectral subspaces
(isotypical components) of $A$. Note that the subalgebra of fixed points of $A$ with respect to $\mathfrak a$ can be strictly
smaller than the spectral subspace corresponding to the trivial corepresentation of $\mathfrak G$ (in the CQG case they are
equal). This creates specific problems that we solve in Sections 4,5 and 6 devoted to the proof of our main result which is parallel
to \cite{CY}, Theorem 6.4 and \cite{Nes1}, Theorem 3.3:

\begin{theorem} \label{main}
Let $\mathfrak G$ be a regular coconnected finite quantum groupoid. Then the following two categories are equivalent:

(i) The category of unital $\mathfrak G$-$C^*$-algebras with unital $\mathfrak G$-equivariant $*$-homomorphisms as morphisms.

(ii) The category of  pairs $(\mathcal M, M)$, where $\mathcal M$ is a left module $C^*$-category over $C^*$-tensor category
${\bf UCorep}(\mathfrak G)$ of unitary corepresentations of $\mathfrak G$ and $M$ is a generator in $\mathcal M$, with
equivalence classes of unitary module functors respecting the prescribed generators as morphisms.
\end{theorem}

This proof divides into three parts. First, given a unital $\mathfrak G$-$C^*$-algebra $A$, we show in Section 4 that the category
$\mathcal D_A$ of finitely generated equivariant $C^*$-correspondences whose morphisms are equivariant maps, is a strict left module
category over ${\bf UCorep}(\mathfrak G)$. The algebra $A$ itself is a generator in $\mathcal D_A$. The idea of such a construction
in the CQG case was proposed in \cite{CY}.

Vice versa, it is shown in \cite{Nes1} that any pair $(\mathcal M, M)$ as above generates so-called weak tensor functor. Using
this functor, we construct in Section 5 an algebra whose $C^*$-completion is a unital $\mathfrak G$-$C^*$-algebra. Finally, we
show in Section 6 that the two above mentioned constructions are mutually inverse which gives the equivalence of the categories
in question.

It was shown in \cite{NY} in the CQG case that ${\bf UCorep}(\mathfrak G)$-module categories parameterized by unitary tensor
(not weak tensor !) functors correspond to Yetter-Drinfeld $\mathfrak G$-$C^*$-algebras. In a subsequent work we expect to get
a similar result for finite quantum groupoids and to apply it to the description of coideal $C^*$-subalgebras of quotient
type.

Our standard references are: \cite{ML} for general categories, \cite{EGNO} for tensor categories, \cite{Nes} for $C^*$- and
$C^*$-tensor categories, \cite{L} for Hilbert $C^*$-modules, and \cite{NV} for finite quantum groupoids.
\end{section}
%%%%%%%%%%%%%%%%%%%%%%%%%%%%%%%%%%%%%%%%%%%%%%%%%%%%%%%%%%%%%%%%%%%%%%%%%%%%%%%%%%%%%%%%%%%%%%%%%%%%%%%%%%%%%%%%%%%%%%%
\begin{section}{Finite quantum groupoids, their  representations, comodules and 
corepresentations}

{\bf 1. Finite quantum groupoids} A {\em weak Hopf  $C^*$-algebra} $\mathfrak
G=(B,\Delta,S,\varepsilon)$ is a finite dimensional $C^*$-algebra $B$ with the
comultiplication $\Delta : B\to B\otimes B$, counit $\varepsilon : B\to C$, and
antipode $S:B\to B$ such that $(B,\Delta, \varepsilon)$ is a coalgebra and the
following axioms hold for all $b,c,d\in B$ :
\begin{enumerate}
\item[(1)] $\Delta$ is a (not necessarily unital)  $*$-homomorphism :
$$
\Delta(bc) = \Delta(b)\Delta(c), \quad \Delta(b^*) =
\Delta(b)^*,
$$
\item[(2)] The unit and counit satisfy the identities (we use the Sweedler leg notation $\Delta(c)=c_1\otimes c_2,\
(\Delta\otimes id_B)\Delta(c)=c_1\otimes c_2\otimes c_3$ etc.):
\begin{eqnarray*}
\varepsilon(bc_1)\varepsilon(c_2d) &=&\varepsilon(bcd), \\
(\Delta(1)\otimes 1)(1\otimes \Delta(1)) &=& (\Delta\otimes id_B)\Delta(1),
\end{eqnarray*}
\item[(3)]
$S$ is an anti-algebra  and anti-coalgebra  map such that
\begin{eqnarray*}
m(id_B \otimes S)\Delta(b) &=& (\varepsilon\otimes id_B)(\Delta(1)(b\otimes 1)),\\
m(S\otimes id_B)\Delta(b) &=& (id_B \otimes \varepsilon)((1\otimes b)\Delta(1)),
\end{eqnarray*}
where $m$ denotes the multiplication.
\end{enumerate}
\medskip

The right hand sides of two last formulas are called {\em target}
and {\em source counital maps} $\varepsilon_t$ and $\varepsilon_s$,
respectively. Their images are unital $C^*$-subalgebras of $B$ called
{\em target} and {\em source counital subalgebras} $B_t$ and $B_s$,
respectively.

The dual vector space $\hat B$ has a natural structure of a weak Hopf
$C^*$-algebra $\hat{\mathfrak G}=(\hat B,\hat\Delta,\hat S,\hat\varepsilon)$
given by dualizing the structure operations of  $B$:
\begin{eqnarray*}
<\varphi\psi,\, b> &=& < \varphi\otimes\psi,\, \Delta(b)>, \\
<\hat\Delta(\varphi),\, b\otimes c> &=& <\varphi,\, bc>, \\
<\hat S(\varphi),\, b> &=& < \varphi,\, S(b)>, \\
< \phi^*,b> &=& \overline{ < \varphi,\, S(b)^*> },
\end{eqnarray*}
for all $b,c\in B$ and $\varphi,\psi\in \hat B$. The unit of $\hat B$ is
$\varepsilon$ and the counit is $1$.

The counital subalgebras commute elementwise, we have $S\circ\varepsilon_s =
\varepsilon_t\circ S$ and $S(B_t) =B_s$. We say that $B$ is {\em connected} if
$B_t \cap Z(B)= \mathbb{C}$ (where $Z(B)$ is the center of $B$), coconnected if
$B_t \cap B_s = \mathbb{C}$, and {\em biconnected} if both conditions are
satisfied.
\medskip

The antipode $S$ is unique, invertible, and satisfies $(S\circ *)^2 =id_B$. We
will only consider {\em regular} quantum groupoids, i.e., such that
$S^2|_{B_t}=id$. In this case, there exists a canonical positive element $H$ in
the center of $B_t$ such that $S^2$ is an inner automorphism implemented by $G=
HS(H)^{-1}$, i.e., $S^2(b) = GbG^{-1}$ for all $b\in B$. The element $G$ is
called the canonical group-like element of $B$, it satisfies the relation\
$\Delta(G) =(G\otimes G)\Delta(1)= \Delta(1)(G\otimes G)$.

There exists a unique positive functional $h$ on $B$,
called a {\em normalized Haar measure} such that
$$
(id_B\otimes h)\Delta = (\varepsilon_t\otimes h)\Delta,\quad
h\circ S =h,\quad h\circ \varepsilon_t = \varepsilon,\quad (id_B\otimes h)\Delta(1_B) = 1_B.
$$
We will dehote by $H_h$ the GNS Hilbert space generated by $B$ and $h$ and by $\Lambda_h:
B\to H_h$ the corresponding GNS map.
\vskip 0.5cm
{\bf 2. Unitary representations}

By definition, the objects of the category $URep(\mathfrak G)$ of unitary representations of $\mathfrak G$
are left $B$-modules of finite rank such that the underlying vector space is a Hilbert space $H$ with a scalar
product $<\cdot,\cdot>$ such that
$$
<b\cdot v,w>=<v,b^*\cdot w>,\quad\text{for all}\quad v,w\in H,\ b\in B,
$$
and morphisms are $B$-linear maps. It is a semisimple linear category whose
simple objects are irreducible $B$-modules. It is also a tensor category: for
objects $H_1,H_2\in URep(\mathfrak G)$, define their tensor product as the
Hilbert subspace $\Delta(1_B)\cdot(H_1\otimes H_2)$ of the usual tensor
product together with the action of $B$ given by $\Delta$. Here we use the   fact that $\Delta(1_B)$ is an orthogonal projection.

The tensor product of morphisms is the restriction of the usual tensor product
of $B$-module morphisms. Let us note that any $H\in URep(\mathfrak G)$ is
automatically a $B_t$-bimodule via $z\cdot v\cdot t:=zS(t)\cdot v,\ \forall
z,t\in B_t, v\in E$, and that the above tensor product is in fact
$\otimes_{B_t}$, moreover the $B_t$-bimodule structure for $H_1\otimes_{B_t} H_2$ is given by $z\cdot \xi \cdot t =(z\otimes S(t))\cdot \xi,\ \forall
z,t\in B_t, \xi \in H_1\otimes_{B_t} H_2$.

One deduces that the above tensor product is associative :
$$
(H_1\otimes_{B_t} H_2)\otimes_{B_t}H_3 =  H_1\otimes_{B_t} (H_2\otimes_{B_t} H_3),
$$
so the associativity isomorphisms are trivial. The unit object of
$URep(\mathfrak G)$ is $B_t$ with the action of $B$ given by $b\cdot z:=
\varepsilon _t(bz),\ \forall b\in B, z\in B_t$ and the scalar product $<z,t>=
h(t^*z)$. The left and right unit morphisms are:
\begin{equation} \label{unit}
l_E(z\otimes_{B_t} v)=z\cdot v\quad\text{and}\quad r_E(v\otimes_{B_t} z)=S(z)\cdot v,\quad\forall z\in B_t, v\in E.
\end{equation}

For any morphism $f:H_1\to H_2$, define $f^*:H_2\to H_1$ as the adjoint linear
map: $<f(v),w>=<v,f^*(w)>,\ \forall v\in H_1, w\in H_2$, it is easy to check
that $f^*$ is $B$-linear. It is clear that $f^{**}=f$, that $(f\otimes_{B_t}
g)^*=f^*\otimes_{B_t} g^*$, and that $End(H)$ is a $C^*$-algebra, for any object
$H$. So $URep(\mathfrak G)$ is a strict finite $C^*$-multitensor category
(i.e., has all the properties of a $C^*$-tensor category except for one: ${\bf
1}$ is not necessarily simple).

In order to make $URep(\mathfrak G)$ a rigid $C^*$-tensor category in the sense
of \cite{Nes}, Definition 2.1.1, we have to define the conjugate for any $H\in
URep(\mathfrak G)$. Take the dual vector space $\hat H$ which is naturally
identified ($v\mapsto \overline v$) with the conjugate Hilbert space $\overline
H: <\overline v,\overline w>=<w,v>,\ \forall v,w\in H$. The action of $B$ on
$\overline H$ is defined by $b\cdot\overline v=\overline
{G^{1/2}S(b)^*G^{-1/2}\cdot v}$, where $G$ is the canonical group-like element
of $\mathfrak G$. Then the rigidity morphisms defined by
\begin{equation} \label{rigid}
R_H(1_B)=\Sigma_i (G^{1/2}\cdot \overline e_i\otimes_{B_t}\cdot e_i),\ \overline R_H(1_B)=\Sigma_i (e_i\otimes_{B_t} G^{-1/2}\cdot\overline e_i),
\end{equation}
where $\{e_i\}_i$ is any orthogonal basis in $H$, satisfy all the needed properties - see \cite{BSz}, 3.6. Also, it is known that the
$B$-module $B_t$ is irreducible if and only if $B_s\cap Z(B)=\mathbb C 1_B$, i.e., if $\mathfrak G$ is connected. So that, we have

\begin{proposition} \label{C*}
$URep(\mathfrak G)$ is a strict rigid finite $C^*$-multitensor category. It is $C^*$-tensor if and only if $\mathfrak G$ is connected.
\end{proposition}
\vskip 0.5cm
{\bf 3. Unitary comodules}

\begin{definition} \label{ucomod}
A  right unitary $\mathfrak G$-comodule is a pair  $(H,\mathfrak a)$, where $H$ is a Hilbert space with scalar product $<\cdot,\cdot>$,
$\mathfrak a : H \to H \otimes B$ is a bounded linear map between Hilbert spaces $H$ and $H\otimes H_h=H\otimes \Lambda_h(B)$, and such that:

(i) $(\mathfrak a \otimes id_B)\mathfrak a= (id_H\otimes \Delta)\mathfrak a$;

(ii) $(id_H \otimes \varepsilon)\mathfrak a= id_H$;

(iii) $<v^1,w>v^2=<v,w^1>S(w^2)^*,\ \ \forall v,w\in H$.

A morphism of unitary $\mathfrak G$-comodules $H_1$ and $H_2$ is a
linear map $T:H_1\to H_2$ such that $\mathfrak a_{H_2}\circ T=(T\otimes
id_B)\mathfrak a_{H_1}$ (i.e., a $B$-colinear map).

Right unitary $\mathfrak G$-comodules with {\bf finite dimensional} underlying Hilbert spaces and their morphisms
form a category which we denote by $UComod(\mathfrak G)$.

We say that two unitary $\mathfrak G$-comodules are equivalent (resp.,
unitarily equivalent) if the space of morphisms between them contains an
invertible (resp., unitary) operator.
\end{definition}
In what follows, we will use the leg notation $\mathfrak a(v)=v^1\otimes v^2$,
for all $v\in H$.
\begin{example} \label{ucoid}
Let us equip a right coideal $I\subset B$ with the scalar product
$<v,w>:=h(w^*v)$. Then the strong invariance of $h$ gives:
$$
<v^1,w>v^2=(h\otimes id_B)((w^*\otimes 1_B)\Delta(v))=
$$
$$
=(h\otimes S^{-1})(\Delta(w^*)(v\otimes 1_B))=<v,w^1>S(w^2)^*.
$$
\end{example}

\begin{remark}
\label{injective}
By (ii) any coaction $\mathfrak a$ is injective
\end{remark}

If $(H,\mathfrak a)$ is a right  unitary $\mathfrak G$-comodule, then $H$ is naturally a unitary
left $\hat{\mathfrak G}$-module via
\begin{equation} \label{mult}
\hat b\cdot v:=v^1<\hat b,v^2>,\ \ \forall \hat b\in\hat B,\ v\in H.
\end{equation}
The unitarity follows from the calculation
$$
<\hat b\cdot v,w> = <v^1<\hat b,v^2>,w> =  <\hat b,<v^1,w>v^2> =
$$
$$
= <\hat b,<v,w^1>S(w^2)^*>=<v,w^1\overline{<\hat b,S(w^2)^*>}>=<v,(\hat b)^*\cdot w>,
$$ for all $v,w\in H$ and $\hat b\in\hat B$. In particular $H$ is a $\hat
B_t$-bimodule.

Due to the canonical identifications $B_t\cong\hat B_s$ and $B_s\cong\hat
B_t$ given by the maps $z \mapsto\hat z=\varepsilon (\cdot z)$ and $t\mapsto
\hat t =\varepsilon (t\cdot)$, $H$ is also a $B_s$-bimodule via $z\cdot v\cdot
t=v^1\varepsilon (z v^2 t)$, for all $z,t\in B_s,\ v\in V$. The maps
$\alpha,\beta:B_s\to B(H)$ defined by $\alpha(z)v:=z\cdot v$ and
$\beta(z)v:=v\cdot z$, for all $z\in B_s, v\in H$ are a $*$-algebra
homomorphism and antihomomorphism, respectively, with commuting images.
Indeed, for instance, for all $v,w\in H, z\in B_s$, one has:
$$
<\alpha(z)v,w>:=<v^1\varepsilon (zv^2),w>=\varepsilon (<v^1,w>zv^2)=
$$
$$
=\varepsilon (<v,w^1>zS(w^2)^*)=<v,w^1>\overline{\varepsilon (S(w^2)z^*)}=
$$
$$
=<v,w^1\varepsilon (S(z^*)w^2)>=<v,\alpha(z^*)w^1\varepsilon (w^2)>=<v,\alpha(z^*)w>.
$$
So that, $\alpha(z)^*=\alpha(z^*)$, and similarly for the map $\beta$.
We have the following useful relations:
\begin{equation} \label{compatible}
\mathfrak a(\alpha(x)\beta(y)v)=v^1\otimes xv^2y\ \forall v\in H, x,y\in B_s.
\end{equation}
and
\begin{equation} \label{compatible'}
\alpha(x)\beta(y)v^1\otimes v^2=v^1\otimes S(x)v^2S(y)\ \forall v\in H, x,y\in B_s.
\end{equation}

The  correspondence (\ref{mult}) is bijective as one has the inverse formula:
if $(b_i)_i$ is a basis for $B$ and $(\hat b_i)$ is its dual basis in $\hat B$,
then set:
\begin{equation}
\label{mult'}
\mathfrak a(v) = \underset{i} \sum (\hat b_i \cdot v) \otimes b_i \ \ \forall v \in H.
\end{equation}
Moreover,  formulas (\ref{mult}) and  (\ref{mult'}) imply also a bijection of
morphisms. Thus, we have two functors, $\mathcal F_1:UComod(\mathfrak G)\to
URep(\hat{\mathfrak G})$ and $\mathcal G_1:URep(\hat{\mathfrak G})\to
UComod(\mathfrak G)$, which are mutually inverse. So, these categories are
isomorphic as linear categories, and we can transport various additional
structures from $URep(\hat{\mathfrak G})$ to $UComod(\mathfrak G)$.

For instance, let us define tensor product of two unitary $\mathfrak
G$-comodules, $(H_1,\mathfrak a_{H_1})$ and $(H_2,\mathfrak a_{H_2})$. As a
vector space, it is
$$
H_1\otimes_{\hat{B_t}} H_2:=\hat\Delta(\hat{1})(H_1\otimes H_2)={\hat{1}}_1\cdot H_1\otimes {\hat{1}}_2\cdot H_2,
$$
and is generated by the elements $ x\otimes_{\hat B_t} y := \hat \Delta (\hat 1)\cdot(x \otimes y)$, where $x\in H_1,
y\in H_2$, so it can be identified with $ H_1\otimes_{{B_s}} H_2$ (see \cite{Pf1}, 2.2 or \cite{Ni}, Chapter 4).
\begin{lemma} \label{complexe}
If $(H_1,\mathfrak a),\ (H_2,\mathfrak b)\in UComod(\mathfrak G)$, then the
projection $P: H_1 \otimes H_2 \to H_1\otimes_{\hat{B_t}} H_2$ defined by $P(v) = \hat \Delta(\hat 1)\cdot v$, for all
$v \in H_1 \otimes H_2$, satisfies
$$
P(x\otimes y) = x^1 \otimes y^1\varepsilon (x^2y^2),\quad\text{for\ all}\ x \in H_1,\ y \in H_2.
$$
\end{lemma}
The proof is the direct calculation using the axiom (2) of a weak Hopf algebra:
$$
{\hat{1}}_1\cdot x\otimes {\hat{1}}_2\cdot y=(x^1\otimes y^1)\varepsilon (x^2 1_1)\varepsilon (1_2 y^2)=
$$
$$
=(x^1\otimes y^1)\varepsilon (x^2y^2).
$$
\begin{corollary} The linear map $\mathfrak a \otimes_{{B_s}} \mathfrak  b$ given by:
$$
v \otimes_{B_s} w\mapsto v^1\otimes_{B_s} w^1\otimes v^2w^2,\ \forall v\in H_1, w\in H_2,
$$
is a coaction of $\mathfrak G$ on $H_1\otimes_{{B_s}} H_2$ (i.e., satisfies Definition \ref{ucomod}, (i), (ii)).
\end{corollary}
\begin{dm} $\forall v\in H_1, w\in H_2$, one has:
\begin{align*} ((\mathfrak a \otimes_{{B_s}}& \mathfrak  b) \otimes i_B)
(\mathfrak a \otimes_{{B_s}} \mathfrak  b)(v \otimes_{B_s} w)=\\
&= ((\mathfrak a \otimes_{{B_s}} \mathfrak  b) \otimes i_B)(v^1\otimes_{B_s} w^1\otimes v^2w^2) \\
&= (\hat \Delta(\hat 1) \otimes 1_{B\otimes B})\cdot(\hat \Delta(\hat 1).(\mathfrak a(v^1)^1\otimes \mathfrak b(w^1)^1)\otimes \mathfrak a(v^1)^2\mathfrak b(w^2) ^2\otimes v^2w^2 )\\
&= (\hat \Delta(\hat 1) \otimes 1_{B\otimes B})\cdot( \mathfrak a(v^1)^1\otimes \mathfrak b(w^1)^1\otimes \mathfrak a(v^1)^2\mathfrak b(w^2) ^2\otimes v^2w^2 ) \\
&= (\hat \Delta(\hat 1) \otimes 1_{B\otimes B})\cdot((\mathfrak a \otimes i_B)\mathfrak a(v))_{134}(\mathfrak b \otimes i_B)\mathfrak b(w))_{2 34})
 \\
&= (\hat \Delta(\hat 1) \otimes 1_{B\otimes B})\cdot((i_E \otimes  \Delta) \mathfrak a(v))_{13}(i_F\otimes  \Delta) \mathfrak b(v))_{23})\\
&= (\hat \Delta(\hat 1) \otimes 1_{B\otimes B}) (i_{EÊ\otimes F} \otimes  \Delta)((\hat \Delta(\hat 1)\otimes 1_B). (\mathfrak a(v))_{13} \mathfrak b(v))_{23})\\
&= (id_{H_1\otimes_{{B_s}} H_2} \otimes  \Delta)(\mathfrak a \otimes_{{B_s}} \mathfrak  b)(v \otimes_{B_s} w)
\end{align*}
Moreover, using Lemma \ref{complexe}, we have :
$$
(id_{H_1\otimes_{{B_s}} H_2} \otimes\varepsilon )(\mathfrak a \otimes_{{B_s}} \mathfrak  b)(v \otimes_{B_s} w)
= v^1 \otimes w^1\varepsilon (v^2w^2) = P(v \otimes w)= v \otimes_{B_s} w
$$
\hfill\end{dm}

The direct calculation shows that the tensor product coaction is unitary. Thus, $UComod(\mathfrak G)$ is a multitensor
category whose associativity morphisms are trivial, the unit object is $(B_s,\Delta|_{B_s})$. It is simple if and only if $\mathfrak G$ is coconnected.
The left and right unit isomorphisms are:
\begin{equation} \label{unit1}
l_H : B_s \otimes_{B_s} H \to H, \ \ z \otimes_{B_s} v \mapsto z\cdot v, \ r_H : H \otimes_{B_s} B_s\to  H, \ \ v \otimes_{B_s} z \mapsto v\cdot z.
\end{equation}
One can check that these isomorphisms are unitary and their inverses are:
\begin{equation} \label{unitinv}
l_H^{-1}(v)=1_1\otimes_{B_s} v^1\varepsilon (1_2v^2)\quad\text{and}\quad r_H^{-1}(v)=v^1\otimes_{B_s}\varepsilon _s(v^2).
\end{equation}

Let us define the conjugate object for $(H,\mathfrak a)\in UComod(\mathfrak G)$. The corresponding Hilbert space is $\overline H$.
In what follows, we use the Sweedler arrows ${\hat b}\rightharpoonup b:=b_1<\hat b,b_2>,\  b\leftharpoonup{\hat b}:=b_2<\hat b,b_1>,
\ \forall b\in B,{\hat b}\in\hat B$.

\begin{lemma} \label{harpoon}
The conjugate object for $(H,\mathfrak a)$ in $UComod(\mathfrak G)$ is $(\overline H,\tilde{\mathfrak a})$, where
$$
\tilde{\mathfrak a}(\overline v)=\overline{v^1}\otimes [{\hat G}^{-1/2}\rightharpoonup (v^2)^*\leftharpoonup{\hat G}^{1/2}],
$$
and $\hat G$ is the canonical group-like element of the dual quantum groupoid $\hat{\mathfrak G}$.
\end{lemma}
\begin{dm}
The unitarity of $\mathcal G_1(\overline H,\tilde{\mathfrak a})$ means that $<\hat b\cdot{\overline v},{\overline w}>_{\overline H}=
<{\overline v},\hat b^*\cdot{\overline w}>_{\overline H}$, for all $v,w\in H$. The left hand side equals to $<{\overline v}^1,
{\overline w}>_{\overline H}<\hat b,{\overline v}^2>$. And the right hand side equals to
$$
<\overline v,\overline{\hat G^{1/2}{\hat S}(\hat b^*)^*\hat G^{-1/2}\cdot w}>_{\overline H}=<\hat G^{1/2}\hat S(\hat b^*)^*\hat G^{-1/2}\cdot w, v>_H=
$$
$$
=<w,\hat G^{-1/2}{\hat S}({\hat b}^*){\hat G}^{1/2}\cdot v>_H=<w,v^1>_H \overline{<\hat G^{-1/2}{\hat S}({\hat b}^*){\hat G}^{1/2},v^2>}=
$$
$$
=<\overline v^1,\overline w>_{\overline H}{<\hat G^{-1/2}\hat b{\hat G}^{1/2},(v^2)^*>}=
$$
$$
=<\overline v^1,\overline w>_{\overline H}<\hat b, [{\hat G}^{-1/2}\rightharpoonup (v^2)^*\leftharpoonup{\hat G}^{1/2}]>.
$$
Comparing the above expressions, we have the result.
\hfill\end{dm}

The rigidity morphisms are given by (\ref{rigid}) with $B_t$ replaced by $B_s$. For any morphism $f$, $f^*$ is the conjugate linear map
on the corresponding Hilbert spaces, the colinearity of $f$ implies that $f^*$ is colinear. So that, we have

\begin{proposition} \label{C**}
$UComod(\mathfrak G)$ is a strict rigid finite $C^*$-multitensor category. It is $C^*$-tensor if and only if $\mathfrak G$ is
coconnected.
\end{proposition}
\vskip 1cm {\bf 4. Unitary corepresentations}
\begin{definition} \label{corepresentation}
A {\bf right unitary corepresentation} of $\mathfrak G$ on a Hilbert space $H$ is a partial isometry $V\in B(H)\otimes B$ such that:

(i) $V_{12}V_{13} =(id_{B(H)}\otimes\Delta)(V)$.

(ii) $(id_{B(H)}\otimes\varepsilon)(V)=id_{B(H)}$.

If $U$ and $V$ are two right corepresentations on Hilbert spaces $H_U$ and $H_V$, respectively, a morphism between them is a
bounded linear map $T\in B(H_U,H_V)$ such that $(T \otimes 1_B)U= V(T\otimes 1_B)$. The vector space of such morphisms is denoted by
$Mor(U,V)$. We will denote by $UCorep(\mathfrak G)$ the category whose objects are right unitary corepresentations $(H,V)$ on {\bf
finite dimensional} vector spaces with morphisms as above.

One  says that $U$ and $V$ are equivalent (resp., unitarily equivalent) if $Mor(U,V)$ contains an invertible (resp., unitary) operator.
\end{definition}
\begin{proposition} \label{comcorep}
If $(H,\mathfrak a)$ is a unitary $\mathfrak G$-comodule, let us define an operator $V$ on
$H\otimes H_h$ as follows:
$$
V(x\otimes \Lambda_h y):=x^1\otimes \Lambda_h(x^2y)),\ \text{for\ all}\ x\in H,\ y\in B.
$$
Then  $V$ is a unitary corepresentation of $\mathfrak G$ on $H$, and one has:
$$
 V^*(x\otimes \Lambda_h y):=x^1\otimes \Lambda_h(S(x^2)y)),\ \text{for\ all}\ x\in H,\ y\in B.
$$
\end{proposition}
\begin{dm}
Let $I_h$ be an implementation of $\Delta$ (for example, $I_h \in B(H_h \otimes H_h):
\Lambda_{h\otimes h}(y' \otimes y) \mapsto \Lambda_{h\otimes h}(\Delta(y)(y' \otimes 1_B))$, see
\cite{Val1}, 3.2) for details), then one has for all $x \in H, y,c \in B$:
\begin{align*}
V_{12}V_{13} (x \otimes \Lambda_h y \otimes &\Lambda_h c)
= (V \otimes 1_B)(x^1\otimes \Lambda_h y \otimes \Lambda_h(x^2c)) \\
&= x^1\otimes \Lambda_h(x^2y) \otimes \Lambda_h(x^3c)\\
&= x^1 \otimes \Lambda_{h\otimes h}(\Delta(x^2)(y\otimes c)) \\
&= x^1 \otimes I_h(x^2\otimes 1_B)I_h^*(\Lambda_{h\otimes h}(y\otimes c))\\
&= (1_{B(H)}\otimes I_h)(x^1\otimes \{(x^2\otimes 1_B)I^*_h\Lambda_{h\otimes h}(y\otimes c)\})\\
&= (1_{B(H)}\otimes I_h)(V \otimes 1_B)(x\otimes I^*_h(\Lambda_{h\otimes h} y \otimes c)) \\
&= (1_{B(H)}\otimes I_h)(V \otimes 1_B)(1_{B(H)}\otimes I_h)^*(x\otimes\Lambda_{h\otimes h} (y \otimes c))\\
&= (1_{B(H)}\otimes \Delta)(V )(x \otimes \Lambda_h y \otimes \Lambda_h c).
\end{align*}
Next, we have, for any decomposition $V = \underset {i \in I} \sum v_i \otimes b_i$ ($v_i \in B(H), b_i \in B)$:
$$
(id_{B(H)}\otimes\varepsilon)(V)(\xi)= (id_{B(H)}\otimes \varepsilon )(\underset {i \in I} \sum v_i(\xi) \otimes b_i) =
$$
$$
=\underset {i \in I} \sum \varepsilon (b_i)v_i(\xi)= (id_{B(H)}\otimes \varepsilon)\mathfrak a(\xi) = \xi,\ \ \forall\xi\in H.
$$
In order to show that $V$ is a partial isometry, consider the separability
element $e_s=(id_B\otimes S)\Delta(1_B)$ of the algebra $B_s$ and the idempotents $e_{\beta,id}=(\beta\otimes id_B)(e_s)
\in\beta(B_s)\otimes B_s$ and $e_{\alpha,S}=(\alpha\otimes S)(e_s)\in\alpha(B_s)\otimes B_s$. As $\alpha$ and $\beta$ are $*$-maps,
these idempotents are orthogonal projections on $H\otimes H_h$. It is straightforward to check, using (\ref{compatible}) and
(\ref{compatible'}), that:
\begin{itemize}
\item for all $x,y \in B_s$, one has:
\begin{equation} \label{comp}
V(\alpha(x)\beta(y)\otimes 1_B)=(1_{B(H)}\otimes x)V(1_{B(H)}\otimes y),
\end{equation}
\begin{equation} \label{comp'}
(\alpha(x)\beta(y)\otimes 1_B)V=(1_{B(H)}\otimes S(x)))V(1_{B(H)}\otimes S(y)).
\end{equation}
\item
$Ve_{\beta,id} = V$, $e_{\alpha,S}V = V$.
\end{itemize}
Moreover, $V$ is invertible in $B(e_{\beta,id}(H\otimes H_h), e_{\alpha,S}(H\otimes H_h))$. Indeed,
consider an operator $W$ acting on $H\otimes H_h$ defined by
$$
W(v\otimes\Lambda_h(b):=v^1\otimes\Lambda_h(S(v^2)b),\ \ \forall v\in H,\ b\in B.
$$
Then we have:
$$
WV(v\otimes\Lambda_h(b):=W(v^1\otimes\Lambda_h(v^2b))=
$$
$$
=v^1\otimes\Lambda_h(S(v^2)v^3b)=v^1\otimes\Lambda_h(\varepsilon _s(v^2)b).
$$
On the other hand,
$$
e_{\beta,id}(v\otimes\Lambda_h(b))=(v\cdot 1_1)\otimes\Lambda_h(S(1_2)b)=v^1\otimes\Lambda_h(S(1_2)\times
$$
$$
\times\varepsilon (v^2 1_1)b)=v^1\otimes\Lambda_h(1_1)\varepsilon (v^2S(1_2))b)=v^1\otimes\Lambda_h(\varepsilon _s(v^2)b).
$$
And similarly $VW=e_{\alpha,S}$, so that $W$ is the inverse of $V$. Finally, we compute, for all $v,w\in H, b,c\in B$:
\begin{align*}
<V(v\otimes\Lambda_h(b)),V(w\otimes\Lambda_h(c))>
&=<v^1\otimes\Lambda_h(v^2b),w^1\otimes\Lambda_h (w^2c)>\\
&=<v^1,w^1> h(c^*(w^2)^*v^2b)\\
&=<v,w^1>h(c^*(w^3)^*[S(w^2)]^*b)\\
&=<v,w^1>h(c^*[\varepsilon _s(w^2)]^*b).
\end{align*}
On the other hand,
$$
<v\otimes\Lambda_h(b),e_{\beta,id}(w\otimes\Lambda_h(c)>=<v\otimes\Lambda_h(b),(w\cdot 1_1)\otimes\Lambda_h(S(1_2)c)>=
$$
$$
=<v,w\cdot 1_1>h(c^*[S(1_2)]^*b)=<v,w^1>\overline{\varepsilon (w^21_1)}h(c^2[S(1_2)]^*b).
$$
These expressions are equal because $\varepsilon _s(x):=1_1\varepsilon (x1_2)=S(1_2)\varepsilon (xS(1_1))=S(1_2)\varepsilon (x1_1))$, for all $x\in B$.
We used above the equality $\varepsilon (xS(z))=\varepsilon (xz)$, for all $z\in B_t$ which can be obtained by applying $\varepsilon \otimes\varepsilon $
to both sides of the equality $\Delta(1_B)(S(z)\otimes 1_B)=\Delta(1_B)(1_B\otimes z)$. As $e_{\beta,id}$ is an orthogonal projection, this means that
$V$ is bounded and $V^*V=e_{\beta,id}$.

Similar reasoning shows that $V^*$ equals to the above mentioned $W$.\hfill\end{dm}

We also have a converse statement.
\begin{proposition} \label{corepcom}
Any unitary corepresentation $V$ of $\mathfrak G$ on a Hilbert space $H$ generates a unitary comodule $(H,\mathfrak a)$, where
$\mathfrak a(v)=V(v\otimes\Lambda_h(1_B))\ \forall v\in H$.
\end{proposition}
\begin{dm}
The first two conditions of Definition \ref{ucomod} follow from  the first two conditions of Definition \ref{corepresentation}.
The relation between $V$ and the coaction $\mathfrak a:v\mapsto v^1\otimes v^2$ is given by $V(v\otimes\Lambda_h(b))=
v^1\otimes\Lambda_h(v^2b)$.
We have seen already that the operator $W$ acting on $H\otimes H_h$ and defined by
$$
W(v\otimes\Lambda_h(b))=v^1\otimes\Lambda_h(S(v^2)b),\ \ \forall v\in H,\ b\in B,
$$
satisfies the relations $VW=e_{a,S}$ and $WV=e_{b,id}$. As $V$ is a partial isometry with initial and final Hilbert
subspaces $e_{a,S}(H\otimes H_h)$ and $e_{b,id}(H\otimes H_h)$, respectively, we have $W=V^*$. Then for all $v,w\in H$ and $ c\in B$,
the equality
$$
<V(v\otimes\Lambda_h(1_B)),w\otimes\Lambda_h(c)>=<v\otimes\Lambda_h(1_B),V^*(w\otimes\Lambda_h(c))>,
$$
can be rewritten as
$$
<v^1,w>h(c^*v^2)=<v,w^1>h(c^*[S(w^2)]^*),
$$
which implies the unitarity of the $\mathfrak G$-comodule in question.
\hfill\end{dm}

Let $(H_1,\mathfrak a)$ and $(H_2,\mathfrak b)$ be two unitary $\mathfrak G$-comodules, and let
$T$ be in $\mathcal B(H_1,H_2)$ intertwining $\mathfrak a$ and $\mathfrak b$, then one has, for all $x \in H_1, b \in B$:

\begin{align*}
V_{H_2}(T \otimes 1)(x \otimes \Lambda_h(b))
&=  (Tx)^1 \otimes \Lambda_h((Tx)^2b) \\
&= (1_{H_2} \otimes \pi'(b))((Tx)^1 \otimes \Lambda_h((Tx)^2) \\
&= (1_{H_2} \otimes \pi'(b))(id_F \otimes \Lambda_h)(\mathfrak  b(Tx)) \\
&=  (1_{H_2} \otimes \pi'(b))(id_F \otimes \Lambda_h)((T \otimes 1)\mathfrak  a(x)) \\
&=  (1_{H_2} \otimes \pi'(b))(T \otimes 1)(id_{H_1} \otimes \Lambda_h)(\mathfrak  a(x))\\
&=  (T \otimes 1)(1_{H_1} \otimes \pi'(b))(id_{H_1} \otimes \Lambda_h)(\mathfrak  a(x))\\
&=  (T \otimes 1)V_{H_1}(x \otimes \Lambda_h(b)).
\end{align*}
Hence, $T \in Mor(V_{H_1},V_{H_2})$.

\begin{corollary}
\label{agit}
The correspondence $\mathcal F_2 $ defined by $\mathcal F_2(H,\mathfrak a)= (V,H)$ and $\mathcal F_2(T) = T$ for all objects
$(H,\mathfrak a)$ and morphisms $T$ of $UComod(\mathfrak G)$, is a functor from $UComod(\mathfrak G)$ to $UCorep(\mathfrak G)$
viewed as semisimple linear categories. The correspondence $\mathcal G_2$ between unitary corepresentations of $\mathfrak G$
and $\mathfrak G$-comodules given by Proposition \ref{corepcom} clearly extends to morphisms and defines a functor inverse to
$\mathcal F_2$, so $UComod(\mathfrak G)$ and $UCorep(\mathfrak G)$ are isomorphic as linear categories. Then we can equip
$UCorep(\mathfrak G)$ with tensor product and duality by transporting these structures from $Comod(\mathfrak G)$.
\end{corollary}
\vskip0.2cm
If $(U,H_U),\ (V,H_V)\in UCorep(\mathfrak G)$, let us define their tensor product.
\begin{lemma}
\label{coin}
One has $(P \otimes id_B)U_{13}V_{23} = U_{13}V_{23}(P \otimes id_B) = U_{13}V_{23},$
where $U_{13}V_{23}\in B(H_U\otimes H_V)\otimes B$ and $P$ was defined in Lemma \ref{complexe}.
\end{lemma}
\begin{dm}
There exist finite families $\{b_k\}$ and $\{b'_k\}$ in $B_s$ such that $\Sigma_k b'_kb_k =\Sigma_k b_kb'_k  = 1_B$,
and for all $x \in H_U$ and all $y \in H_V$ one has:
$$
P(x\otimes y) = \hat \Delta(\hat 1)\cdot(x \otimes y) =\Sigma_k \beta(b_k)x \otimes \alpha'(b'_k)y,
$$
where $\alpha'$ is the $*$-representation of $B_s$ corresponding to $(V,H_2)$. Using four times  (\ref{comp'}), one has:
$$
(P \otimes id_B)U_{13}V_{23}= \Sigma_k (\beta(b_k)  \otimes \alpha'(b'_k)\otimes id_B)U_{13}V_{23}=
$$
$$= \Sigma_k (\beta(b_k)  \otimes id_{H_V}\otimes id_B)U_{13}(id_{H_U} \otimes \alpha'(b'_k)\otimes id_B)V_{23}$$
$$= \Sigma_k (\beta(b_k)  \otimes id_{H_V}\otimes id_B)U_{13}(id_{H_U} \otimes id_{H_V} \otimes S(b'_k)V_{23}$$
$$= \Sigma_k (\beta(b_k)\beta(b'_k) \otimes id_{H_V}\otimes id_B)U_{13}V_{23}$$
$$= \Sigma_k (\beta(b'_k b_k) \otimes id_{H_V}\otimes id_B)U_{13}V_{23} = U_{13}V_{23}$$
$$= \Sigma_k U_{13}(id_{H_U} \otimes id_{H_V} \otimes b_kb'_k)V_{23} $$
$$= \Sigma_k U_{13}(\beta(b_k) \otimes id_{H_V} \otimes 1_B )V_{23}(id_{H_U} \otimes \alpha'(b'_k) \otimes id_B)$$
$$= \Sigma_k U_{13} V_{23}(\beta(b_k) \otimes \alpha'(b'_k) \otimes id_B)= U_{13} V_{23}(P\otimes id_B).$$
\end{dm}

Lemma \ref{coin} justifies the following:

\begin{definition} If $(U,H_U),\ (V,H_V)\in UCorep(\mathfrak G)$, their tensor product is the bounded linear map:
$$
U \otop V = U_{13}V_{23} = (P \otimes id_B)U_{13}V_{23}(P \otimes id_B)
$$
viewed as an element of $B(H_U\otimes_{B_s} H_V)\otimes B$.
\end{definition}

\begin{proposition}
$U \otop V\in UCorep(\mathfrak G)$, it acts on $H_U\otimes_{B_s} H_V$ and:
$$
\mathcal G_2(U,H_1) \otimes_{B_s} \mathcal G_2(V,H_2) = \mathcal G_2(U \otop V,H_1\otimes_{B_s} H_2).
$$
\end{proposition}
\begin{dm}
If $(U,H_U),\ (V,H_V)\in UCorep(\mathfrak G)$, let $U = \underset i \sum u_i \otimes b_i$, $V =  \underset j \sum u_j \otimes b_j$
be decompositions of $U$ and $V$. Then $U \otop V = \underset{i,j} \sum u_i \otimes v_j \otimes b_ib_j$, and let us define $\theta_{U \otop V }
\in B( H_U\otimes_{B_s} H_V, H_U\otimes_{B_s} H_V \otimes B)$ by:
$$
\theta_{U \otop V }( x \otimes_{B_s} y) =  \underset{i,j} \sum u_i \otimes v_j (P(x \otimes y))\otimes b_ib_j.
$$
Then, using Lemma \ref{coin}, one has:
\begin{align*}
\theta_{U \otop V }( x \otimes_{B_s} y)
&=  \underset{i,j} \sum u_i \otimes v_j (x \otimes y))\otimes b_ib_j
 =  \underset{i,j}  \sum P(u_i (x)\otimes v_j (y)) \otimes b_ib_j \\
 &= \mathfrak (a_U \otimes \mathfrak a_V)(x\otimes_{B_s} y),
\end{align*}
and the result follows.\hfill\end{dm}
\vskip 0.5cm

The unit object $\bf 1$ of $ UCorep(\mathfrak G)$ with respect to $\otop$ acts on $B_s$ and is defined by $z\otimes b\mapsto 1_1\otimes 1_2zb$,
for all $z\in B_s,\ b\in B$. It is simple if and only if $\mathfrak G$ is coconnected. The conjugate object for $(V,H)\in UCorep(\mathfrak G)$
is the unitary corepresentation acting on $\overline H$ via $\overline V(\overline x\otimes\Lambda_h(y))=\overline x^1\otimes
\Lambda_h((\overline x^2)^*y)$, where $\tilde{\mathfrak a}(\overline x)$ is described in Lemma \ref{harpoon}, and the rigidity morphisms are
the same as in $UCorep(\mathfrak G)$. For any morphism $f$, again $f^*$ is the conjugate bounded linear map on the corresponding Hilbert spaces.
So that, we have
\begin{proposition} \label{C**1}
$UCorep(\mathfrak G)$ is a strict rigid finite $C^*$-multitensor category. It is $C^*$-tensor if and only if $\mathfrak G$ is coconnected.
\end{proposition}

The simple objects of this category are exactly irreducible corepresentations of $\mathfrak G$. Let us denote by $\Omega$ the set of equivalence
classes of irreducibles and choose a representative $U^x$ in any class $x\in \Omega$. The regular corepresentation of $\mathfrak G$ is decomposed
as follows:

\begin{equation}\label{decomp}
W=\oplus_{x\in \Omega} dim(x) U^x,
\end{equation}
where $dim(x)$ is the dimension of the Hilbert space on which $U^x$ acts.

\begin{definition} \label{matrcoef}
Let $(U,H_U)\in UCorep(\mathfrak G)$ and $\{m_{i,j}\}_{i,j=1}^n$ be the matrix units of $B(H_U)$ with respect to some
orthonormal basis $\{e_i\}_{i=1}^n$ in $H_U$. Then
$$
U=\Sigma_{i,j=1}^n m_{i,j}\otimes U_{i,j},
$$
where $U_{i,j}\ (i,j=1,...,n)$ are called the matrix coefficients of $U$ with respect to $\{e_i\}$. Put $B_U:=Span(U_{i,j})^n_{i,j=1}$;
in particular, we denote $B_{U^x}$ by $B_x$.
\end{definition}

\begin{remark} \label{mcprop}
Let us summarize some properties of  matrix coefficients of $U^x\ (x\in \Omega)$ which can be proved in a standard way.

(i) $B_{\oplus^p_{k=1} U_k}=span\{B_{U_1},...,B_{U_p}\}$ for any finite direct sum of unitary corepresentations. In particular,
(\ref{decomp}) implies that $B=\oplus_{x\in \Omega} B_x$.

(ii) Decomposition $U\otop V=\oplus_z d_z U^z$ with multiplicities $d_z$ implies that $B_UB_V\subset \oplus_z B_z$, where
$z$ parameterizes the irreducibles of the above decomposition.

(iii) The definition of a unitary corepresentation written in terms of $U^x_{i,j}$:
$$
\Delta(U^x_{i,j})= \Sigma^{dim(x)}_{k=1} U^x_{i,k}\otimes U^x_{k,j} ,\quad \varepsilon (U^x_{i,j})=\delta_{i,j},\quad U^x_{i,j}=S(U^x_{j,i})^*,
$$
for all $i,j=1,...,dim(x)$, gives: $B_x\otimes B_x=\Delta(1_B)(B_x\otimes B_x)$, $\Delta(B_x)\subset \Delta(1_B)(B_x\otimes B_x)$ and $B_U=S(B_U)^*$.
We also have $B_{\overline U}=(B_U)^*$.
\end{remark}
\begin{example} \label{beps}
In the case of the trivial corepresentation of $\mathfrak G$ associated with $(\Delta_{\mid B_s}, B_s)$, we will use the notation $B_\varepsilon$
instead of $B_U$.  Let $\{b_i\}^{dim B_s}_{i=1}$ be an orthonormal basis in
$B_s$ with respect to the scalar product $<z,t>=\varepsilon (t^*z)\ \forall z,t\in B_s$. Then one can write $\Delta(1_B)=\Sigma^{dim B_s}_{i=1}
b_i^*\otimes S(b_i)$ (see \cite{NV}, 2.3.3), which implies: $\Delta(b^*_j)=\Sigma^{dim B_s}_{i=1}(b_i^*\otimes S(b_i)b^*_j)$, so
$U^{\varepsilon }_{i,j}=S(b_i)b^*_j$, for all $i,j=1,...,dim B_s$. This means that $B_\varepsilon $ is the unital $C^*$-algebra $B_tB_s$.
\end{example}

{\bf 5. Fiber functor and reconstruction theorem}
\vskip 0.5cm
Let $Q$ and $R$ be two unital $C^*$-algebras. By definition, a $(Q,R)$-correspon- \newline dence is a right Hilbert $R$-module $\mathcal E$ (see \cite{L})
with a unital $*$-homomorphism $\varphi:Q\to\mathcal L(\mathcal E)$, where $\mathcal L(\mathcal E)$ is the $C^*$-algebra of all bounded $R$-linear adjointable
operators on $\mathcal E$. If $Q=R$, we call it an $R$-correspondence. $R$-correspondences form a $C^*$-multitensor category $Corr(R)$ with interior
tensor product $\otimes_R$ and adjointable $R$-bilinear maps as morphisms.

There exists another definition of a $(Q,R)$-correspondence, due to Alain  Connes, this   is a triple $(H, \alpha, \beta)$ where $H$ is a Hilbert space
equipped with unital $*$-homomorphism $\alpha:Q\to B(H)$ and $*$-anti-homomorphism $\beta: R\to B(H)$ whose images commute in $B(H)$. Then $H$ is
a $(Q,R)$-bimodule via $q\cdot v\cdot r:=\alpha(q)\beta(r)v$, for all $q\in Q,r\in R, v\in H$.

In this paper, we are especially interested in the particular case, when $Q=R$ is a finite dimensional $C^*$-algebra equipped with a faithful
tracial state $\phi$. Below we treat this particular case in detail.

\begin{lemma} \label{cor} Both definitions of an $R$-correspondence are equivalent.
\end{lemma}

\begin{dm} (i) If $(H,\alpha,\beta)\in Corr(R)$, define, for any $\eta\in H$, an operator $\Pi(\eta):H_{\phi}\to H$ by $\Pi(\eta)\Lambda_\phi(r):=
\beta(r)\eta$, for all $r\in R$, where $H_\phi$ is the GNS Hilbert space generated by $(R,\phi)$. Then define an $R$-valued scalar product:
$$
<\xi,\eta>_R:=\Pi(\xi)^*\Pi(\eta),\quad\text{for\ all}\quad \xi, \eta\in H.
$$
It is clear that $<\xi,\eta>_{R}$ is in fact in $\pi_\phi(R)$. Finally, $\Pi(\beta(b)\eta)=\Pi(\eta)\pi_\varepsilon (b)$, so $<\xi,\beta(b)\eta>_{R}=<\xi,\eta>_{R}\pi_\phi(b)$, for all $\xi, \eta \in H, b\in R$. Moreover, together with the unital
$*$-representation $\alpha$ we have on $H$ the structure of an $R$-correspondence in the sense of the first definition.

(ii) Vice versa, if $H$ is an $R$-correspondence in this last sense, then one can define a usual scalar product
$<\xi,\eta>=\phi(<\eta,\xi>_R)$, for all $\eta,\xi\in H$, and there are clearly  a unital $*$-homomorphism $\alpha:R\to B(H)$ and a unital
$*$-anti-homomorphism $\beta:R\to B(H)$ whose images commute in $B(H)$. Thus, $(H,\alpha,\beta)$ is an $R$-correspondence in the sense of A. Connes.
\hfill\end{dm}

A morphism between $(H,\alpha,\beta)$ and $(K,\alpha',\beta')$ is a map $T\in B(H,K)$ intertwining $\alpha$ and $\alpha'$ and also $\beta$
and $\beta'$, then $Corr(R)$ is a semisimple linear category. If $(H, \alpha, \beta), (K, \alpha', \beta')\in Corr(R)$, we define their tensor product:
$$
(H, \alpha, \beta)Ê\otimes_R (K, \alpha', \beta') =  ((\beta \otimes \alpha')(e)(H Ê\otimes  K),  \alpha \otimes 1_K, 1_H \otimes \beta'),
$$
where $e$ is the symmetric separability idempotent for $R$, so $e_{\beta,\alpha'} = (\beta \otimes \alpha')(e)$ is an orthogonal projection.
For the sake of simplicity we shall denote $H_1\otimes_R H_2 := e_{\beta,\alpha'}(H Ê\otimes  K)$, and $v\otimes _R w = e_{\beta,\alpha'}(vÊ\otimes w)$,
for all $v  \in H, w \in K$. The  unit object is $R$ with the GNS scalar product defined by $\phi$. The unit isomorphisms are as follows:
$$
l_H(z\otimes_R v):=z\cdot v\quad\text{and}\quad r_H(v\otimes_R z):=v\cdot z,\quad\forall z\in R, v\in H.
$$
They are isometric, for example:
$$
||l_H(z\otimes_R v)||^2:=||z\cdot v||_H^2=\phi(1_R)||z\cdot v)||^2=||1_R\otimes(z\cdot v)||^2=||z\otimes_R v||^2.
$$
The conjugate of a morphism $T:H_1\to H_2$ is just the adjoint operator $T^*:H_2\to H_1$, so $Corr(R)$ is a $C^*$-multitensor category.
We denote by $Corr_f(R)$ its full subcategory with finite dimensional underlying Hilbert spaces. The unit object is simple if and only if
$R$ is a full matrix algebra.

For all objects of the three above categories: $URep(\mathfrak G)$, $UComod(\mathfrak G)$, and $UCorep(\mathfrak G)$, the underlying Hilbert
spaces are $B_s$-correspondences, so each of these categories has a forgetful $C^*$-tensor functor with values in $Corr_f(B_s)$.
\vskip 0.5cm
In order to reformulate in suitable terms the reconstruction theorem of Tannaka-Krein type for finite quantum groupoids proved initially in
\cite{Ha}, \cite{SZ}, recall the construction of the canonical {\bf Hayashi functor} $\mathcal H$.

Let $\mathcal C$ be a rigid finite $C^*$-tensor category and $\Omega=Irr(\mathcal C)$ be an exhaustive set of representatives of equivalence
classes of its simple objects. Let $R$ be the $C^*$-algebra $R=\mathbb C^{\Omega}=\underset {x\in\Omega}\bigoplus\mathbb C p_x$, where $p_x=p_x^*$
are mutually orthogonal idempotents: $p_x p_y=\delta_{x,y}p_x$, for all $x,y\in \Omega$. Then $\mathcal H$ is a functor from $\mathcal C$ to
$Corr_f(R)$ defined by:
$$
\mathcal H(x)=H_x=\underset {y,z\in\Omega}\bigoplus \mathcal C(z,y\otimes x),\quad\text{for\ every}\ x\in \Omega,
$$
where $\mathcal C(x,y)$ is the vector space of morphisms $x\to y$. The $R$-bimodule structure on $H_x$ is given by:
$$
p_y\cdot H_x\cdot p_z=\mathcal C(z,y\otimes x),\quad\text{for\ all}\quad x,y,z\in\Omega.
$$
If $y \in \Omega$ and $f\in\mathcal C(x,y)$, then $\mathcal H(f):H_x\to H_y$ is defined by:
$$
\mathcal H(f)(g)=(id_z\otimes f)\circ g,\quad\text{for\ any}\ z,t\in\Omega\ and\ g\in p_z\cdot H_x\cdot p_t.
$$

The inverse natural isomorphisms $J^{-1}_{x,y}: H_x\otimes H_y\to H_x\underset {R}\otimes H_y$ are:
$$
J^{-1}_{x,y}(v\otimes w)=a_{z,x,y}\circ (v\otimes id_y)\circ w\in p_z\cdot H(x\otimes y)\cdot p_t,
$$
for all $v\in p_z\cdot H_x\cdot p_t, w\in p_t\cdot H_y\cdot p_s, z,s,t\in\Omega$. Here $a_{z,x,y}$ are the associativity isomorphisms
of $\mathcal C$.

We define the scalar product on $H_x$ as follows. If $x,y,z \in \Omega$ and $f,g\in \mathcal C(z,y\otimes x)$, then $g^*\in \mathcal C(y\otimes x,z)$
and $g^*\circ f\in End(z)=\mathbb C$, so one can put $<f,g>_x= g^*\circ f $. The subspaces $\mathcal C(z,y\otimes x)$ are declared to be orthogonal,
so $H_x\in Corr_f(R)$. Dually, $\overline H_x\in Corr_f(R)$ via $z_1\cdot{\overline v}\cdot z_2=\overline{z_2^*\cdot v\cdot z_1^*}$, for all $z_1,z_2\in R,
v\in H_x$. Now one can check that $\mathcal H : \mathcal C \to  Corr_f(R)$ is a unitary tensor functor in the sense of \cite{Nes} 2.1.3.

\begin{theorem} \label{reconstruction}
Let $\mathcal C$ be a rigid finite $C^*$-tensor category and $\Omega=Irr(\mathcal C)$. Let $R$ be the $C^*$-algebra $\mathbb C^{\Omega}$ and
$\mathcal H : \mathcal C \to Corr_f(R)$ be the Hayashi functor. Then the vector space
\begin{equation} \label{algebra}
B =  \underset {x \in \Omega} \bigoplus \overline H_x\otimes H_x,
\end{equation}
has a regular biconnected finite quantum groupoid structure $\mathfrak G$ such that $\mathcal C\cong UCorep(\mathfrak G)$ as $C^*$-tensor categories.
\end{theorem}
\begin{dm}
A rigid finite $C^*$-tensor category $\mathcal C$ is semisimple and spherical, so \cite{Pf1}, Theorems 1.1 and 1.2 claims that $B$ has a structure
of a selfdual regular biconnected semisimple weak Hopf algebra. The algebra of the dual quantum groupoid $\hat{\mathfrak G}$ is (see \cite{SZ}, \cite{M}):
\begin{equation} \label{dualalg}
\hat B =  \underset {x\in\Omega}\bigoplus B(H_x),
\end{equation}
the duality is given, for all $x,y\in\Omega, A\in B(H_y),v,w\in H_x$ by:
$$
<A,\overline w\otimes v>=\delta_{x,y}<Av,w>_x.
$$
$\hat B$ is clearly a $C^*$-algebra with the obvious matrix product and involution, its coproduct is given (see \cite{M}  Theorem 1.3.4) by:
$$
\hat \Delta(\hat b) = \underset{i \in I} \sum (s(r_i ) \otimes t(p_i)) J\hat bJ^{-1},\quad\text{for\ any}\ \hat b \in \hat B,
$$
where  $\underset{i \in I} \sum (r_i \otimes p_i)$ is the symmetric separability element of $R$ hence $\underset{i \in I}\sum (s(r_i )\otimes
t(p_i)) = \hat\Delta(\hat 1)$ is an orthogonal projection in $\hat B\otimes\hat B$; moreover $J = \underset{x,y \in \Omega} \bigoplus
\mathcal H_{x,y}$ is a unitary as a direct sum of unitaries. Then one can easily deduce that $\hat\Delta(\hat b^*) =\hat\Delta(\hat b)^*$,
so both $\hat {\mathfrak G}$ and $\mathfrak G$ are finite quantum groupoids.

The explicit structure of $\mathfrak G$ is given in \cite{Pf1}, Theorems 1.1 and 1.2. If $v,w\in H_x, g,h\in H_y$ and $\{e^{x}_j\}$
is an orthogonal basis in $H_x\ (\forall x,y\in\Omega$), then:
\begin{equation} \label{coproduct}
\Delta(\overline w\otimes v) =  \underset {j}\bigoplus (\overline w\otimes e^{x}_j)_x\otimes(\overline{e^{x}_j}\otimes v)_x,
\end{equation}
\begin{equation} \label{counit}
\varepsilon (\overline w\otimes v) =<v,w>_x,
\end{equation}
\begin{equation} \label{product}
(\overline w\otimes v)_x\cdot (\overline g\otimes h)_y =  (\overline{J^{-1}_{x,y}(w\otimes g)}\otimes J^{-1}_{x,y}(v\otimes h))_{x\otimes y}
\in \overline{H_{x\otimes y}}\otimes H_{x\otimes y},
\end{equation}
\begin{equation} \label{unit3}
1_B = \underset {x\in\Omega}\bigoplus (\rho_x\otimes \rho^{-1}_x)_{\bf 1},
\end{equation}
where $\rho_x$ is the unit constraint attached to $x$, so $\rho^{-1}_x\in p_x\cdot H_{\bf 1}\cdot p_x$ and $\rho_x=\overline{\rho^{-1}_x}$.
In order to define the antipode, consider the natural isomorphisms $\Phi_x:H_x\to\overline H_{x^*}$ and $\Psi_x:\overline H_x\to H_{x^*}$ given by:
$$
\Phi_x=\rho_y(id_y\otimes \overline {ev_x})\circ a_{y,x,x^*}\circ(v\otimes id_{x^*}), \Psi_x=(\overline v\otimes id_{x^*})\circ a^{-1}_{y,x,x^*}
\circ(id_y\otimes coev_x)\circ\rho^{-1}_y.
$$
Here $ev_x$ and $coev_x\ (x\in\Omega)$ are the rigidity morphisms. Then we define:
\begin{equation} \label{antipode}
S(\overline w\otimes v)=[\Phi_x(v)\otimes\Psi_x(\overline w)]_{x^*}.
\end{equation}

Any $H_x$ is a right $B$-comodule via
$$
\mathfrak a_x(v)= \underset {j}\Sigma e^x_j \otimes \overline{e^x_j}\otimes v ,\quad\text{where}\quad v\in H_x,
$$
one checks that it is unitary which gives the equivalence $\mathcal C\cong UCorep(\mathfrak G)$.
\hfill\end{dm}
\end{section}

%%%%%%%%%%%%%%%%%%%%%%%%%%%%%%%%%%%%%%%%%%%%%%%%%%%%%%%%%%%%%%%%%%%%%%%%%%%%%%%%%%%%%%%%%%%%%%%%%%%%%%%%%%%%%%%%%%%%%%%%%%%%%%%%%%%%%%
\begin{section}{Coactions of finite quantum groupoids on unital C*-algebras}

{\bf 1. Canonical implementation of a coaction}

\begin{definition}
\label{sieste}
A  right coaction of a finite quantum groupoid $\mathfrak G$ on a unital $*$-algebra $A$, is a $*$-homomorphism
$\mathfrak a: A \to A \otimes B$  such that:

%1) $(A, \mathfrak a)$ is a right $B$-comodule.

1) $(\mathfrak a \otimes i)\mathfrak a= (id_A\otimes \Delta)\mathfrak a$.

2) $(id_A\otimes\varepsilon )\mathfrak a=id_A.$

3) $\mathfrak a(1_A) \in A \otimes B_t$.

One also says that $(A,\mathfrak a)$ is a $\mathfrak G$-$*$-algebra.
\end{definition}

\begin{remark}
If $A$ is a $C^*$-algebra, then $\mathfrak a$ is automatically continuous, even an isometry by \ref{injective} and
\cite{Ped} 1.5.7.
\end{remark}

\begin{proposition}\label{leonid}
Any right coaction of $\mathfrak G$ on a unital $*$-algebra $A$ is simplifiable: the set $\mathfrak a(A)(1_A \otimes B) =
\{\mathfrak a(a)(1_A \otimes b)\ | \ a \in A , b\in B \}$ generates $\mathfrak a (1_A)(A\otimes B)$ as a vector space.
\end{proposition}
\begin{dm}
Using Sweedler notations (which makes sense here as $B$ is finite dimensional), one has:
\begin{align*}
\mathfrak a (1_A)(a \otimes 1_B)
&= (id_A \otimes \varepsilon  \otimes id_B)(\mathfrak a \otimes id_B)[\mathfrak a (1_A)(a \otimes 1_B)]\\
&= (id_A \otimes \varepsilon  \otimes id_B)[(id_A \otimes \Delta)\mathfrak a(1_A)) \mathfrak a (a) \otimes 1d_B)]\\
&= (id_A \otimes \varepsilon  \otimes id_B)[({1_A }^1 \otimes \Delta({1_A }^2)( a^1 \otimes a^2 \otimes 1_B)]\\
&= (id_A \otimes \varepsilon  \otimes id_B)[({1_A }^1 \otimes \Delta(1_B)({1_A }^2 \otimes 1_B)( a^1 \otimes a^2 \otimes 1_B)]\\
&= (id_A \otimes \varepsilon  \otimes id_B)[(1_A \otimes \Delta(1_B))({1_A }^1a^1 \otimes {1_A }^2a^2 \otimes 1_B)]\\
&= (id_A \otimes \varepsilon  \otimes id_B)[(1_A \otimes \Delta(1_B))(a^1 \otimes a^2 \otimes 1_B)]\\
&= (id_A \otimes \varepsilon _t)\mathfrak a(a).
\end{align*}
Definition 2.1.1 (3) of \cite{NV} gives that:
\begin{align*}
\mathfrak a (1_A)(a \otimes 1_B)
&= (id_A \otimes m)(id_A \otimes id_B \otimes  S)(id_A \otimes \Delta)\mathfrak a(a)\\
&= (id_A \otimes m)(id_A \otimes id_B \otimes  S)(\mathfrak a \otimes id_B)\mathfrak a(a)\\
&= (id_A \otimes m)(\mathfrak a(a^1) \otimes  S(a^2)).
\end{align*}
Finally, the trivial equality: $(id_A \otimes m)(x \otimes y \otimes z) = (x \otimes y)(1_A \otimes z)$ implies:
\begin{align*}
\mathfrak a (1_A)(a \otimes 1_B)
&= \mathfrak a(a^1)( 1 \otimes S(a^2)).
\end{align*}
So $\mathfrak a (1_A)(a \otimes 1_B)$ belongs to the vector space generated by $\mathfrak a(A)(1_A \otimes B)$.
\hfill\end{dm}

Let us introduce the unital $*$-homomorphism $\alpha:B_s\to A: \alpha(x):=x\cdot 1_A$. Equalities (\ref{compatible}) and
(\ref{compatible'}) show that, for all $x\in B_s$ and $a\in A$ :
\begin{equation} \label{algcomp}
\mathfrak a(\alpha(x)a)=(1_A\otimes x)\mathfrak a(a),
\end{equation}
\begin{equation} \label{algcomp'}
(\alpha(x)\otimes 1_B)\mathfrak a(a)=(1_A\otimes S(x))\mathfrak a(a).
\end{equation}

It is helpful to note that
\begin{equation} \label{legs}
\mathfrak a(1_A)=(\alpha\otimes id_B)\Delta(1_B).
\end{equation}
Indeed :
$$
\alpha(1_{1})\otimes 1_{2}:=1_{1}\cdot 1_A\otimes 1_{2}=(id_A\otimes\varepsilon )[(1_A\otimes 1_{1})\mathfrak a(1_A)]\otimes  1_{2}=
$$
$$
=1^{1}_A\otimes(\varepsilon \otimes id_B)\Delta(1^{2}_A)=\mathfrak a(1_A).
$$

\begin{lemma} \label{michelin4} (cf. \cite{Val2} 3.1.5, 3.1.6).
If $(A,\mathfrak a)$ is a $\mathfrak G$-$*$-algebra A, then:

(i) The set $A^{\mathfrak a} = \{a \in A | \mathfrak a(a) = \mathfrak a(1_A) (a \otimes 1_B) \}$ is a unital $*$-subalgebra of $A$
(it is a unital $C^*$-subalgebra of $A$ when $A$ is a $C^*$-algebra) commuting pointwise with $\alpha(B_s)$.

(ii) The map $T^\mathfrak a := (id_A \otimes h)\mathfrak a$ (where $h$ is the normalized Haar measure of $\mathfrak G$) is a conditional
expectation from $A$ to $A^\mathfrak a$; it is faithful when $A$ is a $C^*$-algebra.
\end{lemma}
\begin{dm}
(i) For all $a \in A^\mathfrak a$ and $x \in B_s$, one has:
\begin{align*}
\mathfrak a(a\alpha(x))=\mathfrak a(1_A)(a \otimes 1_B)(1_A \otimes x)\mathfrak a(1_A) = \mathfrak a(\alpha(x)a),
\end{align*}
so $A^\mathfrak a$ commutes pointwise with $\alpha(B_s)$, then it is stable with respect to the multiplication and the $*$-operation
in $A$; moreover if $A$ is a $C^*$-algebra, it is clearly norm closed in $A$, so this is a unital $C^*$-subalgebra of $A$.

(ii) Since $h_{\mid{B_t}}=\varepsilon_{\mid{B_t}}$ (see \cite{NV}, 7.3.2), one has $T^\mathfrak a(1_A) := (id_A \otimes h)\mathfrak a(1_A)=1_A$,
from where, for all $a\in A^\mathfrak a$:
$$
T^{\mathfrak a}(a)=(id_A \otimes h)(\mathfrak a(1_A)(a \otimes 1_B)) = (id_A\otimes h)(\mathfrak a(1_A))a =a.
$$
Now, if $E_t = (id_B \otimes h)\Delta$ is the target Haar conditional expectation of $\mathfrak G$, one has, for all $a\in A$:
\begin{align*}
\mathfrak a (T^{\mathfrak a}&(a))
= \mathfrak a((id_A \otimes h)\mathfrak a (a)) = (id_A \otimes id_B \otimes h)((\mathfrak a \otimes id_B)\mathfrak a (a)) \\
&= (id_A \otimes id_B \otimes h)(id_A \otimes \Delta)\mathfrak a (a)  = (id_A \otimes E_t)\mathfrak a (a) \\
&= (id_A \otimes E_t)(\mathfrak a (1_A)\mathfrak a (a))\\
&= (id_A \otimes E_t)(\mathfrak a(1_A)(a^1 \otimes a^2)) = \mathfrak a(1_A)(a^1 \otimes E_t(a^2)) \\
&= \mathfrak a(1_A)(1_A\otimes E_t(a^2))(a^1\otimes 1_B) =\mathfrak a(1_A)(\beta(S(E_t(a^2)))\otimes 1_B)(a^1 \otimes 1_B) \\
&=  \mathfrak a(1_A)(\beta(S(E_t(a^2)))a^1 \otimes 1_B).
\end{align*}
Using the fact proved above that $(id_A \otimes h)(\mathfrak a(1_A)) = 1_A$, this implies that:
\begin{align*}
(id_A \otimes h)\mathfrak a (T^{\mathfrak a}(a))
&=   (id_A \otimes h)\mathfrak a(1_A)(\beta(S(E_t(a^2)))a^1 \otimes 1_B) \\
&= \beta(S(E_t(a_2)))a^1.
\end{align*}
But since $h\circ E_t=h$, one has also:
\begin{align*}
(id_A \otimes h)\mathfrak a (T^{\mathfrak a}(a))
&= (id_A \otimes h)\mathfrak a(1_A)(a_1\otimes E_t(a_2)) \\
&= (id_A \otimes h)(\mathfrak a(1_A)(a_1\otimes a_2)) \\
&= (id_A \otimes h)(\mathfrak a(1_A)\mathfrak a(a)) = T^\mathfrak a(a).
\end{align*}
One deduces that $T^{\mathfrak a}(a) = \beta(S(E_t(a_2)))a^1$ and
\begin{align*}
\mathfrak a (T^{\mathfrak a}(a))
=  \mathfrak a(1_A)(\beta(S(E_t(a_2)))a_1 \otimes 1_B) = \mathfrak a(1_A)( T^{\mathfrak a}(a)\otimes 1_B).
\end{align*}
This implies that $T^{\mathfrak a}(A) = A^{\mathfrak a}$, moreover, $T^\mathfrak a \circ T^\mathfrak a = T^\mathfrak a$.
Finally, for all $c,d \in A^\mathfrak a$ and $a \in A$, one has:
\begin{align*}
T^\mathfrak a (cad)
&= (id_A \otimes h)\mathfrak a(cad) =( id_A \otimes h)(\mathfrak a(c)\mathfrak a(a)\mathfrak a(d)) \\
&= ( id_A \otimes h)((1_B \otimes c)\mathfrak a(a)(1_B \otimes d)) = cT^\mathfrak a(a)d.
\end{align*}
When $A$ is a $C^*$-algebra, $T^\mathfrak a$ is faithful because $\mathfrak a$ and $h$ are faithful.
\hfill \end{dm}
\begin{definition}
\label{michelin5}
Let $(A,\mathfrak  a)$ be a unital $\mathfrak G$-$*$-algebra, then  unital $*$-subalgebra 
$$
A^{\mathfrak a} = \{a \in A / \mathfrak a(a) = \mathfrak a(1_A) (a \otimes 1_B) \}
$$
is called the subalgebra of invariants (or fixed points) of $(A,\mathfrak a)$.
\end{definition}
\vskip 0.5cm
\begin{proposition} \label{invariance}

Let $(A,\mathfrak a)$ be a unital $\mathfrak G$-$C^*$-algebra and $\phi$ be an element in $A^*$, then
the following assertions are equivalent:

i) for any $a \in A$ one has: $(\phi \otimes i_B)\mathfrak a(a) \in B_s$;

ii) $\phi \circ T^\mathfrak a = \phi$;

iii) there exists a linear form $\omega$ on $A^\mathfrak a$ such that $\phi = \omega \circ T^\mathfrak a$;

iv) for any $x,y \in A$, one has:

$$
(\phi\otimes id_B)(\mathfrak a(x)(y\otimes 1_B))=(\phi\otimes S)((x\otimes 1_B)\mathfrak a(y)).
$$
\end{proposition}

\begin{dm}
Clearly, ii) and iii) are equivalent. If ii) is true and if $E_s = (h \otimes i_B)\Delta$ is the
source Haar conditional expectation of  $\mathfrak G$, then i) is true because, for all $\omega' \in B^*$
and $a \in A$, one has:
\begin{align*}
\omega'((\phi \circ i_B)\mathfrak a(a))
&= (\phi \circ \omega')\mathfrak a(a) = (\phi \circ \omega')(T^\mathfrak a \otimes i_B)\mathfrak a(a)
\\
&= (\phi \circ \omega')((id_A \otimes h)\mathfrak a \otimes id_B)\mathfrak a(a) \\
&= (\phi \circ \omega')(id_A \otimes h \otimes id_B)(\mathfrak a \otimes id_B)\mathfrak a(a) \\
&= (\phi \circ \omega')(id_A \otimes h \otimes id_B)(id_A \otimes \Delta)\mathfrak a(a)  \\
&= (\phi \circ \omega')(id_A \otimes (h \otimes id_B)\Delta)\mathfrak a(a)  \\
&= (\phi \circ \omega')(id_A \otimes E_s)\mathfrak a(a)   \\
&= \omega' (E_s((\phi \circ id_B)\mathfrak a(a))).
\end{align*}

If i) is true, one has:
\begin{align*}
\phi(a)
&= \phi((id_A \otimes \varepsilon )\mathfrak a(a)) = \varepsilon ((\phi \otimes id_B)\mathfrak a(a)) = \varepsilon (E_s(\phi \otimes id_B)\mathfrak a(a))) \\
&= (\phi \otimes \varepsilon )(id_A \otimes E_s)\mathfrak a(a) = (\phi \otimes \varepsilon )(id_A \otimes (h \otimes id_B)\Delta)\mathfrak a(a) \\
&= (\phi \otimes \varepsilon )(id_A \otimes h \otimes id_B)(id_A \otimes \Delta)\mathfrak a(a)\\
& = (\phi \otimes \varepsilon)(id_A \otimes h \otimes id_B)
(\mathfrak a \otimes id_B)\mathfrak a(a) \\
& = (\phi \otimes \varepsilon)(T^\mathfrak a \otimes id_B)\mathfrak a(a) = (\phi \circ T^\mathfrak a)(id_A\otimes \varepsilon )\mathfrak a(a) =
(\phi \circ T^\mathfrak a)(a),
\end{align*}
which is ii), so the three first assertions are equivalent.

Further, if iv) is true, then we have, applying it to $x \in A$ and $y = 1_B$:
$$
(\phi\otimes id_B)\mathfrak a(x) =(\phi\otimes S)((x\otimes 1_B)\mathfrak a(1_A)),
$$
which implies i). Suppose now that i) is true (and so ii) and iii) as well).
First, for all $a \in A, z\in B_t$, the equality (\ref{algcomp'}) gives:
$$
a^1S(z)\otimes a^2=a^1\otimes a^2 z.
$$
Next, the equality $y^1\otimes\varepsilon_t(y^2)=(1_A^1y)\otimes 1_A^2$ (which can be proven directly),
the equality $\varepsilon_t(b)=b_1S(b_2), \forall b\in B$ and assertion i)  give:
\begin{align*}
(\phi\otimes id_B)(\mathfrak a(x)(y\otimes 1_B))
&= \phi(x^1y)x^2=\phi(x^11_A^1y)x^21_A^2=\phi(x^1y^1)x^2\varepsilon_t(y^2)\\
&=\phi(x^1y^1)x^2y^2S(y^3)=\phi((xy^1)^1)(xy^1)^2 S(y^3)\\
&=\phi((xy^1)^1)\varepsilon_s((xy^1)^2)S(y^3).
\end{align*}

Now, using the definition of $\varepsilon_s$ and the equality $\varepsilon_t(bz)=
\varepsilon_t(bS(z))$ which is true for all $b\in B, z\in B_t$, we have:
\begin{align*}
(\phi\otimes id_B)(\mathfrak a(x)(y\otimes 1_B))
&= \phi((xy^1)^1)\varepsilon((xy^1)^2 (1_B)_2)(1_B)_1 S(y^3)\\
&=\phi((xy^1)^1)\varepsilon((xy^1)^2 S((1_B)_2))(1_B)_1 S(y^3)\\
&
=( \phi \otimes \varepsilon)(\mathfrak a(xy^1)(1_A \otimes  S((1_B)_2)))(1_B)_1 S(y^3)
\\
&= \phi((i \otimes \varepsilon)(\mathfrak a(xy^1)(1_A \otimes  S((1_B)_2))))(1_B)_1 S(y^3)
\\
&= \phi((xy^1)\cdot S((1_B)_2))(1_B)_1 S(y^3),
\end{align*}
which equals, due to the relation $(ac)\cdot t=a(c\cdot t), \forall a,c\in A$, to:
\begin{align*}
\phi(x(y^1\cdot S((&1_B)_2)))(1_B)_1 S(y^3)
= \phi(x(i\otimes \varepsilon )(\alpha(y^1) (1 \otimes S((1_B)_2))))(1_B)_1 S(y^3)\\
&= \phi(xy^1\varepsilon(y^2S((1_B)_2)))(1_B)_1 S(y^3)\\
&= (\phi \otimes \varepsilon \otimes S)(xy^1 \otimes y^2S((1_B)_2 \otimes y^3S((1_B)_1)\\
&= (\phi \otimes \varepsilon \otimes S)((x\otimes 1 \otimes 1)(\alpha \otimes i) \alpha(y)(1 \otimes \varsigma(S \otimes S)(\Delta(1_B))))\\
&= (\phi \otimes \varepsilon \otimes S)((x\otimes 1 \otimes 1)(i \otimes \Delta ) \alpha(y)(i \otimes \Delta )(1 \otimes 1_B))\\
&= (\phi  \otimes S)((x\otimes1)(i \otimes(i \otimes  \varepsilon)\Delta )\alpha(y))\\
&= (\phi  \otimes S)((x\otimes1)\alpha(y)).
\end{align*}
\hfill\end{dm}
\begin{corollary} \label{centraliseur}
Let $\mathfrak a(1_A) = 1_A^1 \otimes 1_A^2$ be a decomposition of $\mathfrak a(1_A) $ in Sweedler leg notations, and let $\phi$ be a
positive faithful form on $A$ satisfying the conditions of Proposition \ref{invariance}, then $1_A^1$ is in the centralizer of $\phi$.
\end{corollary}
\begin{dm}
Due to i), one has for all $x \in A$: $(\phi\otimes i)\mathfrak a(x)\in B_t$, so $(\phi\otimes S)\mathfrak a(x)=
(\phi \otimes S^{-1})\mathfrak a(x)$, hence by iv) applied twice:
\begin{align*}
(\phi \otimes i)(\mathfrak a(1_A)(x\otimes 1_B))
&= (\phi\otimes S)\mathfrak a(x) = (\phi\otimes S^{-1})\mathfrak a(x) \\
&= (\phi \otimes i)((x \otimes 1_B)\mathfrak a(1_A)),
\end{align*}
which gives the result.
\hfill\end{dm}
\vskip 0.5cm
\begin{definition}\label{invariante}
A linear form  on $A$ satisfying the conditions of Proposition \ref{invariance} is called an invariant form with respect to $\mathfrak a$.
\end{definition}

\begin{example}\label{frayeur}
The Haar measure $h$ is an invariant faithful form on $B$ with respect to the coaction $\Delta$ of $\mathcal G$ on $B$.
\end{example}

\begin{definition}
If $A^{\mathfrak a} = \mathbb C 1_A$, we say that the coaction $\mathfrak a$ is ergodic.
\end{definition}

\begin{example}\label{trivial}
Let $I\subset B$ be a unital right coideal $C^*$-subalgebra with the coaction $\mathfrak a=\Delta|_I$. Then $I^{\mathfrak a}=I\cap B_t$,
so this coaction is ergodic if and only if $I \cap B_t = \mathbb C 1_B$, i.e., if and only if $I$ is connected.
\end{example}

\begin{remark}\label{europe}
Lemma \ref{invariance} iii) shows that the set of $\mathfrak a$-invariant faithful states on $A$ is not empty. Moreover, if $\mathfrak a$ is ergodic,
then the linear form $h_A$ on $A$ defined by $T^\mathfrak a(x) = h_A(x)1_A (\forall x \in A)$ is the unique $\mathfrak a$-invariant faithful state.
\end{remark}

\begin{definition}
Let $H$ be a Hilbert space and $\mathfrak a$ be a coaction of $\mathfrak G$ on a unital $C^*$-subalgebra $A$ of $B(H)$, then an implementation
of $\mathfrak a$ is a unitary corepresentation $V$ of $\mathfrak G$ on $H$ such that, for all $a\in A$, one has:
$$
\mathfrak a(a) = V(a \otimes 1_B)V^*.
$$
\end{definition}
Let us construct a canonical implementation for any coaction.

\begin{proposition} Let $\mathfrak a$ be a coaction of $\mathfrak G$  on $A$ and $\phi$ a faithful
$\mathfrak a$-invariant state on $A$, then the operator $V$ defined on $H_\phi\otimes H_h$ by
$$
V(a\otimes b):=\mathfrak a(a)(1_A\otimes b),\ \text{for\ all}\ a\in A,b\in B,
$$
is a unitary corepresentation of $\mathfrak G$ implementing $\mathfrak a$.
\end{proposition}
\begin{dm}
For the proof that $V$ is a corepresentation of $\mathfrak G$, see the prof of Proposition \ref{comcorep}.
Then Proposition \ref{invariance} and Corollary \ref{centraliseur} imply:
\begin{align*}
<V(a\otimes b),V(a\otimes b)>
&=(\phi\otimes h)((1_A\otimes b^*)\mathfrak a(a^*a)(1_A\otimes b))\\
&=
h[b^*(\phi\otimes id_b)(\mathfrak a(a^*a)(1_A\otimes 1_B))b)]\\
&=h[b^*(\phi\otimes S)[(a^*a\otimes 1_B)\mathfrak a(1_A)]b)\\
&=(\phi\otimes h)[(a^*a1_A^1\otimes b^*S(1_A^2)b)\\
&=<J_\phi \sigma^\phi_{i/2}(1_A^1)^*J_\phi a\otimes S(1_A^2)b , a \otimes b> \\
&=<(J_\phi (1_A^1)^*J_\phi \otimes S(1_A^2))(a \otimes b) , a \otimes b>,
\end{align*}
for all $a\in A,b\in B$, from where
$$
V^*V =( j_\phi \otimes S)\mathfrak a(1_A).
$$
Here $j_\phi (x):= J_\phi x^* J_\phi$ is the Tomita involution associated with $\phi$. Then $V$ is a partial isometry,
by Proposition \ref{leonid} its image is $\mathfrak a(1_A)(H_\phi \otimes H_h)$, so $VV^* = \mathfrak a(1_A)$.
Put $\beta := j_\phi \circ \alpha$, then by Tomita's theory $\beta$ is a faithful anti-representation of $B_s$
whose image commutes in $B(H_\phi)$ with $Im\ \alpha$.

Now, for any  $x,a\in A$ and $b\in B$, one has: $\mathfrak a(x)V(a\otimes b)=\mathfrak a(x)\mathfrak a(a)(1_A\otimes b)=
\mathfrak a(xa)(1_A\otimes b)=V(ax\otimes b)=V(x\otimes 1)(a\otimes b)$. Hence, $\mathfrak a(x)V=V(x \otimes 1)$, and one
deduces that:
\begin{align*}
\mathfrak a(x)
&= \mathfrak a(x) \mathfrak a(1) = \mathfrak a(x)VV^* = V(x \otimes 1)V^*.
\end{align*}
\hfill \end{dm}

\begin{example} \label{coid}
If $I$ is a right coideal *-subalgebra of $B$ and $\Delta|_I$ is a coaction of $\mathfrak G$ on it, the above formula gives
the unitary corepresentation of $\mathfrak G$ which is a canonical implementation of $\Delta$. In particular, if $I=B$
(resp., $I=B_s$), we have the regular (resp., the trivial) unitary corepresentation of $\mathfrak G$.
\end{example}

{\bf 2. Spectral subspaces of $A$}

For any $(U,H_U)\in UCorep(\mathfrak G)$, $H_U$ is a $\mathfrak G$-comodule via $\delta_U:v\mapsto U(v\otimes 1_B)$.
In terms of the matrix coefficients $U_{i,j}\ (i,j=1,...,n)$ with respect to some orthonormal basis $\{e_i\}_{i=1}^n$
in $H_U$, this means that $\delta_U(e_j)=\Sigma_{i=1}^n e_i\otimes U_{i,j}$.

\begin{definition} \label{specsubspace}
Let $A$ be a unital $\mathfrak G$-$C^*$-algebra $A$. We call the spectral subspace of $A$ corresponding to $(U,H_U)$
the linear span $A_U$ of the images of all $\mathfrak G$-comodule maps $H_U\to A$.
\end{definition}

For instance, if $U$ is the trivial corepresentation which is associated with $(\Delta_{\mid B_s}, B_s)$, so  $H_U=B_s$, we will use the notation
$A_\varepsilon$ instead of $A_U$, and we have $\alpha(B_s)\subset A_\varepsilon $. Indeed, $\alpha: B_s\to A$ is a
$\mathfrak G$-comodule map: $\mathfrak a(\alpha(x))=(1_A\otimes x)\mathfrak a(1_A)=(1_A\otimes x)(\alpha\otimes id_B)\Delta(1_B)$
-see (\ref{legs}).

\begin{proposition} \label{speccomod} (cf. \cite{Boca}, Proposition 13).

One can characterize the spectral subspaces as follows:
$$
A_U:=\{a\in A|\mathfrak a(a)\in \mathfrak a(1_A)(A\otimes B_U)\}.
$$
\end{proposition}
\begin{dm}
(i) Let $R:H_U\to A$ is a $\mathfrak G$-comodule map. Then
$$
\mathfrak a(a)=\mathfrak a(R(v))=\mathfrak a(1_A)(R\otimes id)\delta_U(v)\in \mathfrak a(1_A)(R\otimes id)(H_U\otimes B_U),
$$
where $a=R(v),\ v\in H_U$, and
$$
\mathfrak a(1_A)(R\otimes id)(H_U\otimes B_U)\subset \mathfrak a(1_A)(A\otimes B_U).
$$
(ii) Vice versa, let $a\in A$ be such that $\mathfrak a(a)\in\mathfrak a(1_A)(A\otimes B_U)\subset A\otimes B_U$, so
$\mathfrak a(a)=\Sigma_{i,j}(a_{i, j}\otimes U_{i,j})$.
Then, on the one hand,
$$
(\mathfrak a\otimes id_B)\mathfrak a(a)=\Sigma_{i,j}(\mathfrak a(a_{i, j})\otimes U_{i,j}),
$$
and, on the other hand, using Remark \ref{mcprop} (iii),
$$
(\mathfrak a\otimes id_B)\mathfrak a(a)=\Sigma_{i,j}( a_{i, j}\otimes\Delta(U_{i,j}))=
\Sigma_{i,j,k} (a_{i, j}\otimes U_{i,k}\otimes U_{k,j}),
$$
from where $\mathfrak a(a_{k, j})=\Sigma_i (a_{i, j}\otimes U_{i,k})$, for all $k,j=1,...,dim(H_U)$. But $\mathfrak a(1_A)^2=\mathfrak a(1_A)$,
so in fact $\mathfrak a(a_{k, j})=\mathfrak a(1_A)(\Sigma_i a_{i, j}\otimes U_{i,k})$.
We have $a=\Sigma_{j}a_{j,j}$ because the images of both sides of this equality under $\mathfrak a$ coincide and $\mathfrak a$ is injective. So
it suffices to show that any $a_{j,j}$ is the image of some vector from $H_U$ under some $\mathfrak G$-comodule map to $A$. But the map defined
by $e_k\mapsto a_{k,j}$, for all $j,k=1,...,dim(H_U)$ (where $\{e_k\}^{dim(H_U)}_{k=1}$ is the above orthonormal basis in $H_U$), is clearly
a $\mathfrak G$-comodule map and $a_{j,j}$ is the image of the vector $e_j$.
\hfill \end{dm}

\begin{corollary} \label{specprop}

(i) All $A_U$ are closed.

(ii) $A=\oplus_{x\in \Omega} A_{U^x}$.

(iii) $A_{U^x}A_{U^y}\subset \oplus_z A_{U^z},$ where $z$ runs over the set of all irreducible direct
summands of $U^x\otop U^y$.

(iv) $\mathfrak a(A_U)\subset \mathfrak a(1_A)(A_U\otimes B_U)$ and $A_{\overline U}=(A_U)^*$.

(v) $A_\varepsilon $ is a unital $C^*$-algebra.
\end{corollary}

\begin{dm} (i) $\mathfrak a$ is continuous and $dim(B_U)<\infty$, so all $A_U$ are closed.

(ii) Follows from Remark \ref{mcprop} (i).

(iii) Follows from Remark \ref{mcprop} (ii).

(iv) Remark \ref{mcprop} (iii) implies:
$$
\mathfrak a(a^1)\otimes a^2=a^1\otimes\Delta(a^2)\in A\otimes B_U\otimes B_U,
$$
so $\mathfrak a(a^1)\in A\otimes B_U$. As $\mathfrak a(1_A)$ is an idempotent, we have $\mathfrak a(a^1)
\in\mathfrak a(1_A)(A\otimes B_U)$ which means that $a^1\in A_U$. Then the second statement follows.

(v) Follows from Example \ref{beps}
\hfill\end{dm}

\begin{example} \label{aeps}
Let $(\varepsilon ,B_s)$ be the trivial corepresentation of $\mathfrak G$, so $B_\varepsilon =B_sB_t$ is a unital $C^*$-algebra
(see Example \ref{beps}). The definition of $A_\varepsilon $ shows that it is a unital $C^*$-subalgebra of $A$. It contains a
unital $C^*$-subalgebra $\alpha(B_s)A^{\mathfrak a}$ invariant with respect to $\mathfrak a$.
Indeed, if $z\in B_s, a\in A^{\mathfrak a}$, we have, using (\ref{legs}):
$$
\mathfrak a(\alpha(z)a)=(1_A\otimes z)\alpha(1_A)(a\otimes 1_B)\in \alpha(1_A)(\alpha(B_s)A^{\mathfrak a}\otimes B_sB_t).
$$
We will show that for coconnected finite quantum groupoids $A_\varepsilon =\alpha(B_s)A^{\mathfrak a}$.
\end{example}
\end{section}
%%%%%%%%%%%%%%%%%%%%%%%%%%%%%%%%%%%%%%%%%%%%%%%%%%%%%%%%%%%%%%%%%%%%%%%%%%%%%%%%%%%%%%%%%%%%%%%%%%%%%%%%%%%%%%%%%%%%%%

\begin{section} {From coactions to module categories over  $UCorep(\mathfrak G)$ }

{\bf 1. Equivariant $C^*$-correspondences}

The next definition is parallel to the definitions given in \cite{BS1} and \cite{CY}.

\begin{definition} \label{pluton}.
Given a ${\mathfrak G}-C^*$-algebra $(A,\mathfrak a)$, we call a right Hilbert $A$-module $\mathcal E$ $A$-equivariant if it
is equipped with a map $\mathfrak a_{\mathcal E} : \mathcal E \mapsto \mathcal E \otimes B$ such that:

1) $(\mathfrak a_{\mathcal E} \otimes id_B)\mathfrak a_{\mathcal E}= (id_{\mathcal E}\otimes \Delta)\mathfrak a_{\mathcal E}$;
$(id_{\mathcal E} \otimes \varepsilon)\mathfrak a_{\mathcal E}= id_H$;.

2) $\mathfrak a_{\mathcal E}(\xi\cdot a) = \mathfrak a_{\mathcal E}(\xi)\cdot\mathfrak a( a) $, for all $a\in A, \xi\in\mathcal E$.

3) $<\mathfrak a_{\mathcal E}(\xi),\mathfrak a_{\mathcal E}(\eta)>_{A \otimes B} = \mathfrak a(<\xi,\eta>_A)$, for all $\xi,\eta \in \mathcal E$,
where the exterior product $\mathcal E\otimes B$ \cite{L}, Chapter 4, is considered as a right Hilbert $A\otimes B$-module.

Let $\mathcal D_A$ be the category of {\bf finitely generated} $A$-equivariant Hilbert $A$-modules and  morphisms:  equivariant
$A$-linear maps. These maps are automatically adjointable - see \cite{L}, Chapter 1, so $\mathcal D_A$ is a $C^*$-category.

%An element $\mathcal E\in\mathcal D_A$ is said to be irreducible if $End(\mathcal E)=\mathbb C$.
\end{definition}

\begin{remark} \label{simple}
Condition 1) implies that $\mathcal E$ is canonically a $B_s$-bimodule, given by: $x.\xi.y= \xi^1\varepsilon(x\xi^2y)$  $\forall x,y\in B_s,\forall \xi \in \mathcal E$. So $\mathcal E\otimes B$ is a $B_s\otimes B$-bimodule, where $B$ is a
$B$-bimodule via right and left multiplication. Then one proves using (\ref{legs}) and (\ref{compatible'}) that $\mathfrak
a_{\mathcal E}(\xi)\cdot\mathfrak a(1_A)=\mathfrak a_{\mathcal E}(\xi)$, for all $\xi\in\mathcal E$, and that the vector space $(\mathcal
E\otimes B)\cdot\mathfrak a(1_A) $ is generated by $\mathfrak a_{\mathcal E}(\mathcal E)(1_A\otimes B)$ - see the proof of Proposition \ref{leonid}.
\end{remark}

\begin{lemma} \label{bimod}
Any $\mathcal E\in \mathcal D_A$ satisfies the following conditions:

(i) $(z\cdot\zeta)\cdot a=z\cdot(\zeta\cdot a)$, for all $z\in B_s, a\in A$.

(ii) $<z\cdot\zeta,\eta>_A=<\zeta,z^*\cdot\eta>_A$, for all $z\in B_s, \zeta,\eta\in \mathcal E$.
\end{lemma}

\begin{dm} (i) We have:
$$
z\cdot(\zeta\cdot a)=(\zeta\cdot a)^1\varepsilon (z(\zeta\cdot a)^2)=(id_{\mathcal E}\otimes\varepsilon )[((z\cdot\zeta)^1\otimes (z\cdot\zeta)^2)\cdot\mathfrak a(a)]=(z\cdot\zeta)\cdot a.
$$
(ii) The needed equality is equivalent to
$$
\mathfrak a_{\mathcal E}(z\cdot\zeta,\eta>_A)=\mathfrak a_{\mathcal E}(<\zeta,z^*\cdot\eta>_A),
$$
which is the same as
$$
<\mathfrak a_{\mathcal E}(z\cdot\zeta),\mathfrak a_{\mathcal E}(\eta)>_{A\otimes B}=<\mathfrak a_{\mathcal E}(\zeta),\mathfrak a_{\mathcal E}
(z^*\cdot\eta)>_{A\otimes B}
$$
or
$$
<\zeta^1,\eta^1>_A <z\zeta^2,\eta^2>_B=<\zeta^1,\eta^1>_A <\zeta^2,z^*\eta^2>_B.
$$
As we see, the $A$-valued scalar products coincide and both $B$-valued scalar products are equal to $(\zeta^2)^* z^* \eta^2$ which finishes the proof.
\hfill\end{dm}

This lemma shows that any $\mathcal E\in \mathcal D_A$ is automatically a $(B_s, A)$-correspon- dence (see the definition in Section 2);
we call such an object an equivariant $(B_s, A)$-correspondence and denote it by $_{B_s}\mathcal E_{A}$.

%We will show that $\mathcal D_{B_s}\cong UCorep(\mathfrak G)$.

\begin{example} \label{generator}
A ${\mathfrak G}-C^*$-algebra $(A,\mathfrak a)$ itself with the $A$-valued scalar product $<a,b>_A=a^*b\ (\forall a,b\in A)$, is an equivariant
$(B_s,A)$-correspondence.
\end{example}

\begin{theorem} \label{corephilb} If $(V,H_V)$ is a unitary corepresentation of $\mathfrak  G$, then $H_V$ is
an equivariant $B_s$-correspondence ($B_s$ is equipped with the coaction $\Delta|_{B_s}$ of $\mathfrak  G$).
\end{theorem}

\begin{dm} Proposition \ref{corepcom} shows that $(H_V,\mathfrak a_V)$ is a unitary $\mathfrak G$-comodule (where
$\mathfrak a_V(\eta)=V(\eta\otimes\Lambda_h(1_B)),\ \forall \eta\in H_V$) so $H_V$ is a $B_s$-correspondence in the
sense of A. Connes. Then the Hilbert $B_s$-module structure on $H_V$ is described in the proof of Lemma \ref{cor}.

Applying the relations (\ref{compatible}) and (\ref{compatible'}), one has:
$$
\mathfrak a_V(\eta)\cdot \Delta(1_B)= \mathfrak a_V(\eta)\cdot(1_1\otimes 1_2) = \mathfrak a_V(\eta)\cdot (1_B\otimes S(1_1)1_2)
= \mathfrak a_V(\eta),
$$
which implies, for all $\eta\in H_V, t\in B_s$:
$$
\mathfrak a_V(\eta\cdot t) = \mathfrak a_V(\eta)\cdot(\Delta(1_B)(1 \otimes t) )= \mathfrak a_V(\eta)\cdot\Delta(t).
$$
Now, consider $V$ as an element of $B(H_V \otimes H_h)$, where $H_h$ is the GNS Hilbert space constructed by $(B,h)$, the canonical
multiplicative isometry  $I_h$ of $\mathfrak G$ (see \cite{Val2}, Proposition 2.2.4) and its normalized fixed vector $e$ (see \cite{Val1},
\cite{Val5} 2.3 and 2.4)). Applying \cite{Val5}, Lemma 2.1.1,  one has, for all $b' \in B'$ (the commutant of $B$ in $B(H_h)$), $\xi,
\eta \in \mathfrak H$, and $x,x'\in B_s$:
\begin{align*} \label{tenir}
 <\Delta(<\xi,\eta>_\alpha)&(\Lambda_\varepsilon  x \otimes e), \Lambda_\varepsilon  x' \otimes b'e>
 \\
 &= <\Delta(1_B)(1_B \otimes <\xi,\eta>_{B_s})(\Lambda_\varepsilon  x \otimes e), \Lambda_\varepsilon  x' \otimes b'e> \\
 &= (h \otimes \omega_e) ((x'^* \otimes b'^* )\Delta(1_B)(1_B \otimes <\xi,\eta>_{B_s})(x \otimes 1_B)) \\
 &= (h \otimes \omega_e) (\Delta(1_B)(1_B \otimes <\xi,\eta>_\alpha)(xx'^*\otimes b'^*))\\
 &= \omega_e((h \otimes id_B )(\Delta(1_B)(xx'^*\otimes 1_B) <\xi,\eta>_{B_s} b'^*))\\
 &= \omega_e( S(xx'^*) <\xi,\eta>_{B_s} b'^*).
 \end{align*}
On the other hand, taking two decompositions: $V(\xi\otimes e)=\underset {j \in J}\sum (\xi_j \otimes b_je)$ and $V(\eta \otimes e) = \underset {i \in I} \sum(\eta_i \otimes b_ie)$, one computes:
\begin{align*}
&<<\mathfrak a_V(\xi),\mathfrak a_V(\eta)>_{B_s\otimes B}(\Lambda_\varepsilon  x \otimes e), \Lambda_\varepsilon  x' \otimes b'e> = \\
&= \underset {i \in I, j \in J} \sum <<\xi_j \otimes b_j, \eta_i \otimes b_i>(\Lambda_\varepsilon  x \otimes e), \Lambda_\varepsilon  x' \otimes b'e> \\
&= \underset {i \in I, j \in J} \sum <(R(\xi_j)^*R(\eta_i ) \otimes b_j^*b_i)(\Lambda_\varepsilon  x \otimes e), \Lambda_\varepsilon  x' \otimes b'e> \\
&= \underset {i \in I, j \in J} \sum <(R(\xi_j)^*R(\eta_i )\Lambda_\varepsilon  x,\Lambda_\varepsilon  x'><  b_j^*b_i e, b'e> \\
&= \underset {i \in I, j \in J} \sum <R(\eta_i )\Lambda_\varepsilon  x,R(\xi_j)\Lambda_\varepsilon  x'><  b_i e, b_jb'e> \\
&= \underset {i \in I, j \in J} \sum < \beta(x)\eta_i ),\beta(x')\xi_j ><  b_i e, b'b_je> \\
&= < \underset{i\in I}\sum(\beta(x)\otimes 1_B)(\eta_i\otimes b_i e),\underset{j\in J}\sum(\beta(x')\otimes b')(\xi_j \otimes b_j e)> \\
&= <(\beta(x)\otimes 1_B)(V(\eta \otimes e),(\beta(x')\otimes b')V(\xi \otimes e)> \\
&= <V(\eta \otimes S(x)e), V(\xi \otimes S(x')b' e)> \\
& = <e_{\beta,i}(\eta \otimes S(x)e), e_{\beta,i}(\xi \otimes S(x')b'e)> \\
%&=  <\eta \   _{a \circ S}\underset \tau {\otimes\ _i }  S(x)n_0e,   \xi \  _{a \circ S}\underset \tau {\otimes\ _i }S(x')b'e> \\
&=  <<\xi,\eta>_{B_s} S(x)e, S(x')b'e> = <<\xi,\eta>_{B_s} S(xx'^*)e, b'e> \\
&= \omega_e(S(xx'^*)<\xi,\eta>_{B_s} b'^*) = \omega_e(S(xx'^*)<\xi,\eta>_{B_s} b'^*).
\end{align*}
Thus, $<\mathfrak a_V(\xi),\mathfrak a_V(\eta)>_{B_s\otimes B} =\Delta(<\xi,\eta>_{B_s}).$ \hfill \end{dm}

\begin{proposition} \label{recip}

Given an equivariant $B_s$-correspondence $_{B_s}\mathcal E_{B_s}$, define on $\mathcal E$ the scalar product inherited from its $B_s$-scalar product:
$<\xi ,\eta> = \varepsilon (<\eta,\xi>_{B_s})$, for all $\xi,\eta\in\mathcal  E$. Then $V\in B(\mathcal E \otimes H_h)$ defined by
$$
V(\eta \otimes \Lambda_h(b)) =(id_{\mathcal E}\otimes \Lambda_h)(\mathfrak a_{\mathcal E}(\eta)\cdot(1\otimes b)),\ \text{for\ all}\ \eta\in\mathcal E, b \in B,
$$
is a unitary corepresentation of $\mathfrak G$.
\end{proposition}
\begin{dm}
As $\mathcal E$ satisfies the condition 1) of Definition \ref{pluton}, it has a $B_s$-bimodule structure defined by the maps $\alpha, \beta : B_s \to \mathcal L(\mathcal E)$. In particular,
$\beta(n)\xi = \xi\cdot n$, for all $n \in B_s$ and $\xi \in \mathcal E$. Definition \ref{pluton} 2) shows that the right $ B_s$-module structure given by
$\beta$ is the same as the initial $ B_s$-bimodule structure on $\mathcal E$. With the new scalar product on $\mathcal E$, one has:
\begin{align*}
<\beta(n)\xi,\xi>
&= \varepsilon (<\xi,\beta(n)\xi>_{B_s}) = \varepsilon (<\xi,\xi\cdot n>_{B_s}) \\
&= \varepsilon ( <\xi,\xi>_{B_s}n)= \varepsilon (n<\xi,\xi>_{B_s})) \\
&= \varepsilon (<\xi\cdot n^*,\xi>_{B_s})) \\
&= \varepsilon (<\beta(n^*)\xi,\xi>_{B_s})) = <\xi,\beta(n^*)\xi>.
\end{align*}
Hence, $\beta$  is a unital $*$-anti-representation of $B_s$ on $\mathcal E$, and $e_{\beta,i}$ is an orthogonal projection. Moreover, as $\mathcal E$
satisfies the condition 1) of Definition \ref{pluton}, then $V$ defined above satisfies the conditions (i) and (ii) of Definition \ref{corepresentation}
- see the proof of Proposition \ref{comcorep}. On the other hand:
\begin{align*}
<V^*V(\eta \otimes e),& (\eta \otimes e)> = <V(\eta \otimes e), V(\eta \otimes e)> \\
& = <\underset{ i \in I} \sum \eta_i \otimes b_i e, \underset{ i \in I} \sum \eta_i \otimes b_i e> \\
&=(\varepsilon  \otimes h)(<\mathfrak a_{\mathcal E}(\eta),\mathfrak a_{\mathcal E}(\eta)>_{B_s \otimes B})\\
&= (\varepsilon  \otimes h)(\mathfrak a_{triv}(<\eta,\eta>_{B_s })) = h(<\eta,\eta>_{B_s})\\
&= <<\eta,\eta>_{B_s }e,e> = <e_{\beta,i}(\eta \otimes e), \eta\otimes e>.
\end{align*}

As $e$ is separating for $B$, this implies that $V$ is a partial isometry whose initial support is $e_{\beta,i}$.
\hfill \end{dm}

Theorem \ref{corephilb} and Proposition \ref{recip} allow to define two functors : $\mathcal F_3: UCorep(\mathfrak G)\to\mathcal D_{B_s}$ and
$\mathcal G_3:\mathcal D_{B_s}\to UCorep(\mathfrak G)$ on the level of objects, and the morphisms in both cases are just $B$-comodule maps.
These functors are inverse to one another. Indeed, since the $B$-comodule structure is the same in both cases, the only thing to explain is the
relation between the usual scalar product in $H_V$ and the corresponding $B_s$-valued scalar product, but this explanation was done in the proof
of Lemma \ref{cor}. Thus, we have :

\begin{theorem} \label{enfin}
The categories $UCorep(\mathfrak  G)$ and ${\mathcal D}_{B_s}$ are isomorphic.
\end{theorem}

In particular, the unit object ${\bf 1}\in{\mathcal D}_{B_s}$ is $(B_s,\Delta|_{B_s})$ with the $B_s$-valued scalar product
$<b,c>=b^*c$, for all $b,c\in B_s$, and the tensor product is the interior tensor product of $B_s$-correspondences.

{\bf 2. Module categories over $UCorep(\mathfrak G)$ associated with equivariant $C^*$-correspondences.}

\begin{definition} \label{mcat} \cite{CY} Let $\mathcal C$ be a $C^*$-multitensor category with unit object ${\bf 1}$. A $C^*$-category
$\mathcal M$ is called a left $\mathcal C$-module $C^*$-category if there is a bilinear $*$-functor $\boxtimes:\mathcal C\times\mathcal
M\to\mathcal M$ with natural unitary transformations $(X\otimes Y)\boxtimes M\to X\boxtimes (Y\boxtimes M)$ and ${\bf 1}\boxtimes M\to M\
(X,Y\in \mathcal C, M\in \mathcal M)$ making $\mathcal M$ a left module category over $\mathcal C$ - see \cite{EGNO}, Chapter 7. If $\mathcal C$
is strict, we say that $\mathcal M$ is strict (resp., indecomposable) if these natural transformations are identities (resp., if, for all
non-zero $M,N\in \mathcal M$, there is $X\in\mathcal C$ such that $\mathcal M(X\boxtimes M,N)\neq 0$).

We say that an object $M\in\mathcal M$ generates $\mathcal M$ if any object of $\mathcal M$ is isomorphic to a subobject of $X\boxtimes M$
for some $X\in\mathcal C$. $\mathcal M$ is said to be semisimple if the underlying $C^*$-category is semisimple.
\end{definition}
We will always consider $C^*$-categories closed with respect to subobjects, i.e., such that for any object $M$ and any projection
$p\in End(M)$, there are an object $N$ and isometry $v\in \mathcal M(N,M)$ satisfying $p=vv^*$ (if necessary, one can complete given
$C^*$-category with respect to subobjects).

One naturally defines a morphism $F:\mathcal M_1\to\mathcal M_2$ between two $\mathcal C$-module $C^*$-categories as a morphism of the
underlying $C^*$-categories equipped with a unitary natural equivalence $F(X\boxtimes M)\to X\boxtimes F(M),\ \forall\ X\in\mathcal C,\
 M\in\mathcal M$ satisfying some coherence conditions (see \cite{CY}, 2.17).

\begin{lemma} \label{modcat} $\mathcal D_A$ is a strict left module category over $UCorep(\mathfrak G)$ defined
by interior tensor product of $C^*$correspondences over $B_s$.
\end{lemma}
\begin{dm} Given $H_V\in \mathcal D_{B_s}$ and $_{B_s}\mathcal E_{A}\in\mathcal D_A$, equip the vector space $H_V\otimes_{B_s}
\mathcal E$ with $A$-valued scalar product - see \cite{L}, Proposition 4.5:
\begin{equation} \label{tensor}
<v\otimes_{B_s}\zeta,w\otimes_{B_s}\eta>_A=<\zeta,<v,w>_{B_s}\cdot\eta>_A,\ \forall v,w\in H_V,\ \zeta,\eta\in\mathcal E,
\end{equation}
which gives it the $B_s - A$-correspondence structure, and also with the algebraic structure of tensor product of the corresponding
$B$-comodules. One can check that we obtain a new object $H_V\otimes_{B_s}\mathcal E\in\mathcal D_{A}$, and that this construction
is natural both in $V$ and $\mathcal E$. Thus, we have defined a functor $\boxtimes:UCorep(\mathfrak G)\times \mathcal D_A\to
\mathcal D_A$ having the needed properties. Indeed, the first of them is true because $H_{U\otop V}=H_U\otimes_{B_s}H_V$ and because
of the associativity of $\otimes_{B_s}$, and the second one can be proved by direct computation.

Finally, $\boxtimes$ sends adjoint morphisms to adjoint, so it is a $*$-functor.
\hfill\end{dm}

Let us show that $A$ viewed as an object of $\mathcal D_A$ (see Example \ref{generator}) is a generator for $\mathcal D_A$. More precisely,
if $V\in UCorep(\mathfrak G)$, then $H_V\otimes_{B_s} A\in \mathcal D_A$ and the corresponding right coaction of $B$ on $H_V\otimes_{B_s} A$
defines a left action of $\hat B$ on it: $\hat b\cdot v:=v^1<\hat b,v^2>$, for all $v\in H_V\otimes_{B_s} A, \hat b\in \hat B$. If $p\in
\mathcal L(H_V\otimes_{B_s} A)$ is a $\hat B$-invariant orthogonal projection, then one can check that $H_{V,p}=p(H_V\otimes_{B_s} A)$ is a
subobject of $H_V\otimes_{B_s} A$ in $\mathcal D_A$.

\begin{lemma} \label{gen} (cf. \cite{NT}, Lemma 3.2).
For any $\mathcal E\in \mathcal D_A$, there is $V\in UCorep(\mathfrak G)$ and a $\hat B$-invariant projection $p\in
\mathcal L(H_V\otimes_{B_s} A)$ such that $\mathcal E$ is isomorphic to $H_{V,p}$.
\end{lemma}
\begin{dm}
For any fixed $\zeta\in\hat B\cdot\mathcal E=\mathcal E$, the finite dimensional vector space $\hat B\cdot\zeta$ is a $\hat B$-module, so
there is a finite dimensional $\hat B$-submodule $\mathcal E_0$ of $\mathcal E$ such that $\mathcal E_0\cdot A=\mathcal E$. In particular,
there are unital $*$-representations of $B_s\cong\hat B_t$ and $B_t\cong\hat B_s$ on $\mathcal E_0$, so it is a $B_s$-bimodule. Constructing
on this space a $B_s$-valued scalar
product like in the proof of Lemma \ref{cor}, we turn $\mathcal E_0$ into an equivariant $B_s$ correspondence, and Proposition \ref{recip}
allows to construct  $V\in UCorep(\mathfrak G)$ such that the left $\hat B$-modules $H_V$ and $\mathcal E_0$ are isomorphic. Fix an isomorphism
$T_0:H_V\to\mathcal E_0$ and define $T:H_V\otimes_{B_s} A\to\mathcal E$ by $T(v\otimes_{B_s} a)=(T_0 v)\cdot a$. This is a surjective morphism
of $A$-modules. Since $H_V\otimes_{B_s} A$ is a finitely generated Hilbert $A$-module, it makes sense to consider the polar decomposition $T^*=
u|T^*|$. Then $|T^*|$ is an invertible endomorphism of the $A$-module $\mathcal E$, and $u:\mathcal E\to H_V\otimes_{B_s} A$ is an $A$-module
mapping such that $u^*u=\iota$. Property (iii) in Definition \ref{pluton} and non-degeneracy ensure that $T^*$, $|T^*|$, $u=T^*|T^*|^{-1}$, and
$u^*$ are morphisms of $A$-equivariant Hilbert modules. In particular, $u:\mathcal E\to H_{V,p}$ is an isomorphism such that $p=uu^*$.
\hfill\end{dm}

\begin{remark} \label{End} $End_{{\mathcal D}_A}(A)=A^\mathfrak a$. In particular, a coaction $\mathfrak a$ is ergodic if and only if the
generator $A$ of the module category $\mathcal D_A$ is simple.

Indeed, if $T\in End_{{\mathcal D}_A}(A)$, then $\mathfrak a(T(1_A))=(T\otimes id_B)\mathfrak a(1_A)\in A\otimes B_t$. So $T(1_A)\in
A^{\mathfrak a}$ because $(id_A\otimes h)\mathfrak a(T(1_A))=(id_A\otimes\varepsilon )\mathfrak a(T(1_A))=T(1_A)$.

Vice versa, arbitrary $a\in A^{\mathfrak a}$ generates an equivariant endomorphism of $A$ via $T:1_A\mapsto a$.
\end{remark}

We can summarize the above considerations as follows:

\begin{theorem} \label{TK1} Given a regular coconnected finite quantum groupoid $\mathfrak G$, consider two categories:

(i) The category $\mathfrak G-Alg$ of unital $\mathfrak G-C^*$-algebras together with unital $\mathfrak G$-equivariant $*$-homomorphisms as morphisms.

(ii) The category $UCorep(\mathfrak G)-Mod$ of pairs $(\mathcal M,M)$, where $\mathcal M$ is a left $UCorep(G)$-module $C^*$-category and $M$
is its generator, with equivalence classes of unitary $Rep(G)$-module functors respecting the generators as morphisms.

Let us associate with any $\mathfrak G-C^*$-algebra $(A,\mathfrak a)$ the $C^*$-category $\mathcal D_A$ of finitely generated $A$-equivariant
$(B_s,A)$-correspondences with its generator $A$, and with any morphism $f:A_0\to A_1$ in $\mathfrak G-Alg$ the morphism $\mathcal E
\mapsto\mathcal E\otimes_{A_0} A_1$ from $\mathcal D_{A_0}$ to $\mathcal D_{A_1}$. This defines a functor $\mathcal T: \mathfrak G-Alg\to
UCorep(\mathfrak G)-Mod$.
\end{theorem}

The only thing to check is that $\mathcal T$ is well defined on the level of morphisms. This is straightforward because $A_1$ is a left 
$A_0$-module via morphism $f$. This construction was discussed in \cite{CY}, Chapter 7 as "extension of scalars".

%{\bf If $\mathfrak a$ is ergodic, is $\mathcal D_{\mathfrak a}$ semisimple; irreducible as a module category - cf. \cite{CY} ? Vice versa ?}
\end{section}
%%%%%%%%%%%%%%%%%%%%%%%%%%%%%%%%%%%%%%%%%%%%%%%%%%%%%%%%%%%%%%%%%%%%%%%%%%%%%%%%%%%%%%%%%%%%%%%%%%%%%%%%%%%%%%%%%%%%%%%%%%%%%%%%%%%%%%%%%%%%%%%
\begin{section} {From  module categories over $UCorep(\mathfrak G)$ to coactions}

In Sections 5 and 6 we use the approach proposed in \cite{Nes1} with certain modifications reflecting the difference between CQG
and finite quantum groupoids and the fact that we are considering left module categories and right coactions instead of right module
categories and left coactions as in \cite{Nes1}.

\begin{definition} \label{weak} Let $R$ be  a $C^*$-algebra  and let $(\mathcal C,\otimes,{\bf 1})$ be a strict $C^*$-tensor category,
a weak tensor functor from $\mathcal C$ to $Corr(R)$ is a linear functor $F:\mathcal C\to Corr(R)$ together with natural $R$-bilinear
isometries $J=J_{U,V}: F(U)\otimes_{R} F(V)\to F(U\otimes V)$ satisfying the following conditions:

(i) $F({\bf 1})=R$;

(ii) $F(T)^*=F(T^*)$ for any morphism $T$ in $\mathcal C$;

(iii) $J: R\otimes_{R} F(U)\to F({\bf 1}\otimes U)=F(U)$ maps $r\otimes X$ into $Xr$, and
      $J: F(U)\otimes_{R} R\to F(U\otimes{\bf 1})=F(U)$ maps $X\otimes r$ into $rX$, for all $X\in F(U)$;

(iv) $J(id\otimes J)=J(J\otimes id)$;

(v) for all $U,V\in\mathcal C$ and every vector $Y\in F(U)$, the right $R$-linear map $S_Y=S_{Y,U}:F(U)\to F(U\otimes V)$ mapping
$X\in F(U)$ into $J(X\otimes Y)$ is adjointable, and $J(id\otimes S^*_Y)=S^*_Y\circ J$.
\end{definition}

\begin{remark} \label{tensfunctor}
(i) Any unitary tensor functor $F:\mathcal C\to Corr(R)$ is a weak tensor functor - if the conditions (i) - (iv) are satisfied and the maps
$J$ are surjective, then the condition (v) is also satisfied.

(ii) If we consider $F$ as a functor into the category of vector spaces, then $S_Y$ is a natural transformation from $F$ to $F(\cdot\otimes V)$,
and we have
\begin{equation} \label{2.1}
S^*_Y F(T\otimes id)=F(T)\circ S^*_Y,\quad\text{for\ all\ morphisms\ in}\ \mathcal C.
\end{equation}
\end{remark}
We will also need the following modification of \cite{Nes1}, Proposition 3.1:

\begin{proposition} \label{modcatweak}
Let $\mathcal M$ be a strict left module $C^*$-category over a strict  $C^*$-tensor category $\mathcal C$, $M$ be an object in $\mathcal M$,
and denote by $R$ the unital $C^*$-algebra $End(M)$. Then the map $F(U)=\mathcal M(M, U\boxtimes M)\ \forall U\in\mathcal C$ defines a weak tensor
functor $F:\mathcal C\to Corr(R)$, where $X=F(U)$ is a right $R$-module via the composition of morphisms, a left $R$-module via $rX=(id\otimes r)X$,
the $R$-valued inner product is given by $<X,Y>=X^*Y$, the action of $F$ on morphisms is defined by $F(T)X=(T\otimes id)X$, and $J_{X,Y}(X\otimes Y)=(id\otimes Y)X$, for all $X\in F(U), Y\in F(V), X,Y\in\mathcal C$.
\end{proposition}
Let us note that $S_Y(X)=(id\otimes Y)X$ and $S^*_Y(Z)=(id\otimes Y^*)Z$, where $Z\in F(U\otimes V)$.

Now we will describe step by step the reconstruction procedure. Let $\mathcal M$ be a strict left $UCorep(\mathfrak G)$-module $C^*$-category
with generator $M$.

Let $\Omega$ be an exhaustive set of representatives of the equivalence classes of  irreducible objects in $UCorep(\mathfrak G)$.
Consider the following vector space:
\begin{equation} \label{bspace}
A=\underset{x\in\Omega}\bigoplus A_{U^x}:=\underset{x \in\Omega}\bigoplus(F(U^x)\otimes\overline{H_x}),
\end{equation}
and also a much larger vector space:
\begin{equation} \label{recalgebra2}
\tilde{A}=\underset{U\in \|UCorep(G)\|}\bigoplus A_{U}:=\underset{U\in \|UCorep(G)\|}\bigoplus(F(U)\otimes\overline{H_U}),
\end{equation}
where $F(U)=\underset{i}\bigoplus F(U_i)$ corresponds to the decomposition $U=\bigoplus U_i$ into
irreducibles, and $ \|UCorep(G)\| $ is an exhaustive set of representatives of the equivalence classes of objects in  $UCorep(G)$ (these
classes constitute a countable set). $\tilde A$ is a unital associative algebra with the product
$$
(X\otimes\overline\xi)(Y\otimes\overline\eta)=(id\otimes Y)X\otimes(\overline\xi\otimes_{B_s}\overline\eta),
\ \forall (X\otimes\overline\xi)\in A_U, (Y\otimes\overline\eta)\in A_V,
$$
and the unit
$$
1_{\tilde A}=id_M\otimes\overline{1_B}.
$$
Note that $(id\otimes Y)X=J_{X,Y}(X\otimes Y)\in F(U\otop V)$. Then, for any $U\in UCorep(G)$, choose isometries
$w_i:H_i\to H_U$ defining the decomposition of $U$ into irreducibles, and define the projection $p:\tilde A\to A$ by
\begin{equation} \label{p}
p(X\otimes \xi)=\Sigma_i(F(w_i^*) X\otimes\overline{w_i^*\xi}),\ \forall (X\otimes\overline\xi)\in A_U,
\end{equation}
which does not depend on the choice of $w_i$. Indeed, for any other choice of isometries $v_j$ there exists a unitary matrix $u_{ij}$
such that $w_i=\Sigma_{i,j}u_{ij}v_j$. Note also that if $w:H_U\to H_V$ is an isometry between $U,V\in Corep(\mathfrak G)$, then
\begin{equation} \label{pw}
p(F(w) X\otimes\overline{w\xi})=p(X\otimes\overline\xi),\quad \forall (X\otimes\overline\xi)\in A_U.
\end{equation}

\begin{lemma} \label{prod1}
$A$ is a unital associative algebra with the product $x\cdot y:=p(xy)$, for all $x,y\in  A$.
\end{lemma}
\begin{dm}
It suffices to check that $p(p(a)p(b))=p(a)p(b)$, for all $a,b\in\tilde A$. Let $a=(X\otimes\overline\xi)\in A_U,
b=(Y\otimes\overline\eta)\in A_V$, where $U,V\in UCorep(G)$. Choose isometries $u_i$ and $v_j$ corresponding to the decompositions
$U=\bigoplus U_i$ and $V=\bigoplus U_j$ into irreducibles, and let $w_{i,j,k}$ be isometries corresponding to the decomposition of
$U_i\otop V_j$ into irreducibles. Then:
$$
p(a)p(b):=\Sigma_i(F(u_i^*) X\otimes\overline{u_i^*\xi}))\Sigma_j(F(v_j^*) Y\otimes\overline{v_i^*\eta})=
$$
\begin{equation} \label{inter}
=\Sigma_{i,j}((id\otimes F(v_j^*) Y)F(u^*_i) X\otimes\overline{u_i^*\xi\otimes v^*_j\eta})=
\end{equation}
$$
=\Sigma_{i,j,k}(F(w_{i,j,k}^*)(id\otimes F(v_j^*) Y)F(u^*_i) X\otimes\overline{w^*_{i,j,k}(u_i^*\xi\otimes v^*_j\eta)}).
$$
On the other hand, if we apply $p$ to (\ref{inter}), we get the same result.
\hfill\end{dm}

In particular, the vector subspace $A_{\varepsilon }=R\otimes\overline H_{\varepsilon }$ (where $R=End(M)$ and  $H_{\varepsilon}=B_s$) is a unital
$C^*$-subalgebra of $A$ and any $F(U)$ is an $R$-correspondence (see Proposition \ref{modcatweak}).

\begin{lemma} \label{bullet}
If $X$ is in $ F(U)$, then $X^\bullet=S^*_XF(R_U)(1_B)$ is the unique element from $F(\overline U)$ satisfying
$$
<X^\bullet,Y>=F( R^*_U)J(Y\otimes X)\quad\text{for\ all}\ Y\in F(\overline U),
$$
where $R_U$ and $\overline R_U$ come from (\ref{rigid}). We also have:
$$
<X,Y>=F(\overline R^*_U) J(Y\otimes X^\bullet),\quad\forall Y\in F(U).
$$
\end{lemma}
\begin{dm} We compute:
\begin{align*}
<X^\bullet,Y>
&=<S^*_XF(R_U)(1_B),Y>=<F(R_U)(1_B),S_X(Y)>\\
&=F(R^*_U) J(Y\otimes X).
\end{align*}
The uniqueness follows from the faithfulness of the inner product.
As for the last statement, we compute:
\begin{align*}
F(\overline R^*_U) J(Y\otimes X^\bullet)
&=F(\overline R^*_U) J(Y\otimes S^*_XF(R_U)(1_B))\\
&=F(\overline R^*_U)S^*_X J(Y\otimes F(R_U)(1_B))\\
&=S^*_X F(\overline R^*_U\otimes id)F(id\otimes R_U)Y,
\end{align*}
where we have used (\ref{2.1}). The latest expression equals to $S^*_X Y$, where $S^*_X:R\to F(U)$ is given by
$r\to J(X\otimes r)=r\cdot X$, so $S^*_X Y=<X,Y>$.
\hfill\end{dm}
Similarly, for any $\xi\in H_U$ define $\xi^\bullet\in H_{\overline U}$ by
$$
\xi^\bullet=(\overline\xi\otimes id_U)\overline R_U(1_B)=\overline{\hat G^{1/2}\cdot\xi}\ \text{(see (\ref{rigid})),
\ so}\ <\eta,\xi^\bullet>=\overline R_U^*(\xi\otimes\eta)\ \forall \eta\in H_{\overline U},
$$
and consider the map $\bullet:\tilde{A}\to\tilde{A}$
$$
(X\otimes\overline \xi)^{\bullet}:=X^\bullet\otimes\overline{\xi^\bullet}.
$$
\begin{lemma} \label{prod} $A$ is a unital $*$-algebra with the above product and the involution
$x^*:=p(x^{\bullet})$, for all $x\in A$.
\end{lemma}
\begin{dm} First, we prove that $p(p(a)^\bullet)=p(a^\bullet)$, for all $a\in\tilde A$. Take $a=(X\otimes\overline \xi)\in A_U$ and
choose isometries $u_i$ corresponding to the decompositions of $U=\bigoplus U_i$ and into irreducibles. Then
for the standard duality morphisms we have $R_U=\Sigma_i (\overline w_i\otimes w_i)R_i$ and $\overline R_U=\Sigma_i ( w_i\otimes
\overline w_i)\overline R_i$, where $R_i:=R_{U_i}, \overline R_i:=\overline R_{U_i}$. Then
$$
F(R^*_U)(Y\otimes X)=\Sigma_i F(R^*_i)J(F(\overline w^*_i)Y\otimes F(w^*_i)X)=
$$
$$
=\Sigma_i <(F(w^*_i)X)^\bullet,F(\overline w^*_i)Y>,
$$
so $X^\bullet=\Sigma_i F(\overline w_i)(F(w^*_i)X)^\bullet$. Similarly, $\xi^\bullet=\Sigma_i \overline w_i(w^*_i\xi)^\bullet$,
therefore, applying $p$ to $a^\bullet$ and using (\ref{pw}), we have:
$$
p(a^\bullet)=\Sigma_i p(F(\overline w_i)(F(w^*_i)X)^\bullet\otimes\overline{\overline w_i(w^*_i\xi)^\bullet})=
$$
$$
=\Sigma_i p((F(w^*_i)X)^\bullet\otimes \overline{(w^*_i\xi)^\bullet}).
$$
On the other hand, the last expression equals to $p(p(a)^\bullet)$.

Next, in order to prove that $(p(a)p(b))^*=p(b)^*p(a)^*$,
it suffices to prove that $p((a\cdot b)^\bullet)=p(b^\bullet\cdot a^\bullet)$, for all $a,b\in\tilde A$. Take $a=(X\otimes\xi)\in
A_U$ and $b=(Y\otimes\eta)\in A_V$. The unitary $\sigma: H_{\overline V}\otimes
H_{\overline U}\to H_{\overline{U\odot V}}$ mapping $\overline\theta\otimes\overline\zeta$ into $\overline{\zeta\otimes\theta}$
defines an equivalence between $\overline V\otop\overline U$ and $\overline{U\otop V}$, and we have
$$
R_{U\otop V}=(\sigma\otimes id\otimes id)(id\otimes R_U\otimes id)R_V\ \text{and}\ \overline R_{U\otop V}=(id\otimes id\otimes\sigma)
(id\otimes\overline R_V\otimes id)\overline R_U.
$$
Then we compute using Lemma \ref{bullet}, relations $S_{J(X\otimes Y)}=S_Y S_X$ and (\ref{2.1}):
$$
J(X\otimes Y)^\bullet=S^*_{J(X\otimes Y)}F(R_{U\otop V})(1_B)=
$$
$$
=S^*_X S^*_Y F(\sigma\otimes id\otimes id)F(id\otimes R_U\otimes id)F(R_V)(1_B)=
$$
$$
=F(\sigma)S^*_X F(id\otimes R_U)S^*_Y F(R_V)(1_B)=F(\sigma)S^*_X F(id\otimes R_U)(Y^\bullet)=
$$
$$
=F(\sigma)S^*_X J(Y^\bullet\otimes F(R_U)(1_B))=F(\sigma) J(Y^\bullet\otimes S^*_XF(R_U)(1_B))=
$$
$$
=F(\sigma)J(Y^\bullet\otimes X^\bullet).
$$
Similarly, $(\xi\otimes\eta)^\bullet=\sigma(\eta^\bullet\otimes\xi^\bullet)$, from where
$$
(a\cdot b)^\bullet=(F(\sigma)\otimes\overline\sigma)(J(Y^\bullet\otimes X^\bullet)\otimes(\eta^\bullet\otimes\xi^\bullet))
=(F(\sigma)\otimes\overline\sigma)(b^\bullet\cdot a^\bullet).
$$
Applying now $p$, we get $p((a\cdot b)^\bullet)=p(b^\bullet\cdot a^\bullet)$.

In order to show that $**=id$ on A, we will show that $p(a^{\bullet\bullet})=p(a)$, for all $a\in\tilde A$. Take
$a=(X\otimes\xi)\in A_U$ and consider the unitary $u:H_U\to H_{\overline{\overline U}}: \xi\mapsto\overline
{\overline\xi}$. Then $\overline R_{U}=(u\otimes id)\overline R_U$, hence, applying twice Lemma \ref{bullet}, we have: 
$$
<X^{\bullet\bullet},Y>=F(R^*_{\overline U})J(Y\otimes X^\bullet)=F(\overline R^*_{U})F(u^*\otimes id)J(Y\otimes X^\bullet)=
$$
$$
=F(\overline R^*_{U})J(F(u^*)Y\otimes X^\bullet)=<X,F(u^*)Y>,\quad\text{for\ any}\ Y\in F(\overline{\overline U}).
$$
So $X^{\bullet\bullet}=F(u)X$. We also have $\xi^{\bullet\bullet}=\overline{\overline\xi}=u\xi$, from where
$a^{\bullet\bullet}=(F(u)\otimes\overline u)a$, and applying $p$ to both sides of this equality we get $p(a^{\bullet\bullet})=p(a)$.
\hfill\end{dm}

Now define a linear map $\mathfrak a:A\to A\otimes B$ by $\mathfrak a(X\otimes\overline\xi)=X\otimes(-\otimes id_B)U^x({\xi}
\otimes 1_B)$ or, in other words, by
\begin{equation} \label{coact}
\mathfrak a(X\otimes\overline{\xi_i})=X\otimes\Sigma_j(\overline{\xi_j}\otimes U^x_{j,i}),
\end{equation}
where $X\in F(U^x),\ \{\xi_j\}$ is any orthonormal basis in $H_x$ and $U^x_{i,j}$ are matrix coefficients of $U_x$ with respect to
this basis (see Definition \ref{matrcoef}).

\begin{lemma} \label{co}

(i) The map $\mathfrak a$ is a right coaction of $\mathfrak G$ on $A$.

(ii) $A$ admits a unique $C^*$-completion $\overline A$ such that $\mathfrak a$ extends to a continuous coaction
of $\mathfrak G$ on it.
\end{lemma}
\begin{dm} (i) Clearly, $(A,\mathfrak a)$ is a right $B$-comodule. In order to show that $\mathfrak a$ is an algebra homomorphism,
remark that $\tilde A$ is a right $B$-comodule via extension $\tilde{\mathfrak a}$ of $\mathfrak a$ which is defined as in
(\ref{coact}), but with arbitrary $U\in UCorep(\mathfrak G)$. It follows from (\ref{coact}) that $p:\tilde A\to A$ is a
comodule map, and from the formula $U\otop V=U_{13}V_{23}$ that $\tilde{\mathfrak a}$ is a homomorphism, hence $\mathfrak a$
is also a homomorphism.

In order to check that $\mathfrak a$ is $*$-preserving, it suffices to show that $\tilde{\mathfrak a}(a)^{\bullet\otimes *}=\tilde
{\mathfrak a}(a^\bullet)$, for all $a=(X\otimes\overline\xi)\in A_U,\ U\in UCorep(\mathfrak G)$. This is equivalent
to
$$
U(\xi\otimes 1_B)^{\bullet\otimes *}=\overline U(\overline{\hat G^{1/2}\cdot\xi}\otimes 1_B),\quad\forall\xi\in H_U,
$$
which follows from Lemma \ref{harpoon} and a few relations that are easy to check: $(\hat b\cdot v^1)\otimes v^2=
v^1\otimes(v^2\leftharpoonup\hat b)$, $(\hat b\cdot v)^1\otimes(\hat b\cdot v)^2=v^1\otimes(\hat b\rightharpoonup v^2)$,
$(\hat b\rightharpoonup b)^*=\hat S(\hat b)^*\rightharpoonup b^*$, and $(b\leftharpoonup\hat b)^*=b^*\leftharpoonup\hat S(\hat b)^*$,
 for all $v\in H_U,\ b\in B$, and $\hat b\in\hat B$.

Finally, $\mathfrak a(1_A)=id_{\bf 1}\otimes({\overline .}\otimes id_B)U^\varepsilon (1_B\otimes 1_B)=id_{\bf 1}\otimes(-\otimes id_B)
\Delta(1_B)$,  so $\mathfrak a(1_A)\in id_{\bf 1}\otimes\overline{B_s}\otimes B_t$.

(ii) By Lemma \ref{michelin4}, the set $A^{\mathfrak a}$ of all fixed points is a unital $*$-subalgebra of $A$ commuting with $\alpha(B_s)$. 
Moreover, the conditional expectation $T^{\mathfrak a}:=(id_A\otimes h)\mathfrak a$ (where $h$ is the normalized Haar measure of $\mathfrak G$)    
from $A$ onto $A^{\mathfrak a}$ gives rise to a $A^{\mathfrak a}$-valued (pre)inner product for $A$ defined by:
$$
<a,b>_T=T^\alpha(a^*b)\quad\text{for\ all}\quad a,b\in A.
$$
Note that if $a=(X\otimes\overline\xi)\in A_U$, then $T^{\mathfrak a}(p(a))=\Sigma_i(F(w^*_i)X\times\overline{w^*_i\xi})$, where $w_i$ are isometries
corresponding to the decomposition of $U$ into irreducibles such that $\Sigma_i(w_i w^*_i)$ is the projection onto the component of $\varepsilon $.
This implies the mutual orthogonality of the spaces $A_{U^x}\ \forall x\in \Omega$, but $1_A \in A_{\varepsilon}$ hence $T^{\mathfrak a}(A_{U^x}) = 0$ for all $x \not= \varepsilon$. The component
$A_\varepsilon =End(M)\otimes B_s$ is a unital $C^*$-algebra and using (\ref{coact}), by restriction $\mathfrak a$ is a coaction of $\mathcal G$ on the $C^*$-algebra $A_\varepsilon$ and $T^{\mathfrak a}(A_\varepsilon) \subset A_{\varepsilon}$, which implies that  $A^{\mathfrak a} = T^{\mathfrak a} (A) = T^{\mathfrak a} (A_\varepsilon) \subset A_\varepsilon$ and by Lemma \ref{michelin4}, $A^{ \mathfrak a }$ is a unital  $C^*$-subalgebra of $A_\varepsilon$. Therefore  $A$ is a right pre-Hilbert $A^{\mathfrak a}$-module.

The map $T^\alpha$ is completely positive, the $C^*$-algebra $A^{\mathfrak a}$ is unital, and the number of the components $A_{U^x}$ is
finite, so the multiplication on the left gives a faithful $*$-representation $A\to \mathcal L(A)$. One can extend $\mathfrak a$ to the
$C^*$-completion $\overline A$ of $A$ using the reasoning from the proof of \cite{CY}, Proposition 4.4. The map $V$ on $A\otimes B$ defined
by $X(a\otimes b)=\mathfrak a(a)(1_A\otimes b)$, extends (due to the invariance of $h$) to a partial isometry on the right Hilbert
$A^{\mathfrak a}$-module $A\otimes H_h$. The direct calculation shows that the formula $\overline{\mathfrak a}:a\mapsto V(a\otimes
1_B)V^*$ gives the needed extension of the coaction. \hfill\end{dm}

\end{section}
%%%%%%%%%%%%%%%%%%%%%%%%%%%%%%%%%%%%%%%%%%%%%%%%%%%%%%%%%%%%%%%%%%%%%%%%%%%%%%%%%%%%%%%%%%%%%%%%%%%%%%%%%%%%%%%%%%%%%%%%%%%%%%%%%%%%%%%%%%%%%%ù
\begin{section} {Equivalence of categories}

\begin{definition} \label{sspecfun}
Let $(A,\mathfrak a)$ be a unital $\mathfrak G$-$C^*$-algebra and $A_\varepsilon $ be its spectral $C^*$-subalgebra corresponding to the trivial 
corepresentation $\varepsilon $. The {\bf spectral functor} associated with $(A,\mathfrak a)$ is a functor $F:UCorep(\mathfrak G)\to Corr(A_\varepsilon )$ 
defined as follows: for any $U\in UCorep(\mathfrak G)$, put $F(U)=\{X\in H_U\otimes_{B_s} A|U_{13}X_{12}=
(id_A\otimes\mathfrak a)(X)\}=\{X=\Sigma_i(\xi_i\otimes_{B_s} a_i)|\mathfrak a(a_i)=\Sigma_j(a_j\otimes U_{ij}),\ \forall\ i\}$, where $\{\xi_i\}$ is an 
orthonormal basis in $H_U$. Then $F(\varepsilon )=A_\varepsilon $, all $F(U)$ are $A_\varepsilon $-bimodules, and $A_\varepsilon $-valued inner product of $X=\Sigma_i(\xi_i\otimes_{B_s} a_i),\ Y=\Sigma_i(\xi_i\otimes_{B_s} b_i)\in F(U)$ defined by $<X,Y>:=\Sigma_i(a^*_i b_i)$, does not depend on the choice of $\{\xi_i\}$. Putting also $F(T):=T\otimes id$ for morphisms, we have a unitary functor respecting tensor products: if $X=\Sigma_i(\xi_i\otimes_{B_s} a_i)\in F(U),\ Y=\Sigma_j(\eta_j\otimes_{B_s} b_j)\in F(V), \ U,V\in UCorep(\mathfrak G)$, then the maps $J_{U,V}:X\otimes Y\mapsto Y_{23}X_{13}$ are $A_\varepsilon $-bilinear isometries between $F(U)\otimes_{A_\varepsilon }F(V)$ and $F(U\otop V)$.
\end{definition}

\begin{remark} \label{difference}

1) The spectral functor $(F,J)$ associated with a $\mathfrak G$-$C^*$-algebra $(A,\mathfrak a)$ is a weak unitary tensor functor. Indeed, properties (i) - (iv) are immediate, and (v) follows by observing that the adjoint of the map
$$
S_Y:F(U)\to F(U\otop V):\quad X\to Y_{23}X_{13},
$$
is given by $S^*_Y(Z)=Y^*_{23}Z$. Namely, if $Y=\Sigma_i(\eta_j\otimes_{B_s} a_j)$ and $Z=\Sigma_{i,j}(\xi_i\otimes_{B_s} \eta_j\otimes_{B_s} z_{i,j})$ 
for some orthonormal bases $\{\xi_i\}\in H_U$ and $\{\eta_j\}\in H_V$, then
\begin{equation} \label{S*}
S^*_Y Z=\Sigma_{i,j}(\xi_i\otimes_{B_s} a^*_j z_{i,j})\in F(U).
\end{equation}

2) The spectral subspaces $A_U$ can be recovered from $F(U)$ using the canonical surjective maps
$$
F(U)\otimes\overline H_U\to A_U,
$$
which are isomorphisms for irreducible $U$.
\end{remark}

\begin{theorem} \label{spec-weak}
Fix a regular coconnected finite quantum groupoid $\mathfrak G$ and a $C^*$-algebra $C$. By associating to a $\mathfrak G$-$C^*$-algebra 
$(A,\mathfrak a)$ its spectral functor, we get a bijection between isomorphism classes of triples $(A,\mathfrak a,\psi)$, where $\psi:C\to A$ is an 
embedding such that $A_\varepsilon =\psi(C)$, and natural unitary monoidal isomorphism classes of weak tensor functors $UCorep(\mathfrak G)\to 
Corr(C)$.
\end{theorem}
\begin{dm} Isomorphic $\mathfrak G$-$C^*$-algebras produce naturally unitarily monoidally isomorphic weak unitary tensor 
functors, and vice versa. It remains to show that up to some isomorphisms these constructions are mutually inverse.

Let $(A,\mathfrak a)$ be a $\mathfrak G$-$C^*$-algebra with its spectral $C^*$-subalgebra $A_\varepsilon $ corresponding to the trivial corepresentation of $\mathfrak G$, and let $F$ be the associated spectral functor. As $F$ is a weak unitary tensor functor, Lemmas \ref{prod} and \ref{co} allow to construct a unital $\mathfrak G$-$*$-algebra $(A_F,{\mathfrak a}_F)$. One can check that linear
maps sending $p(X\otimes\overline\xi)$ to $(\overline\xi\otimes id)X\in A_U$, for any $(X\otimes\overline\xi)\in F(U)\otimes
\overline H_U\ (U\in UCorep(\mathfrak G))$, define a unital $\mathfrak G$-equivariant homomorphism of algebras. In order to show
that it is $*$-preserving, fix irreducibles $U^x$ and an orthonormal basis $\{\xi_i\}$ in $H_x\ \forall x\in\hat G$. For an
element $X=\Sigma_i(\xi_i\otimes_{B_s} a_i)\in F(U^x)$, we compute, using Lemma \ref{bullet} and identity (\ref{S*}):
$$
X^\bullet=S^*_XF(R_{U^x})(1_B)=S^*_X(\Sigma_j(\overline{\hat G^{-1/2}\xi_j}\otimes\xi_j\otimes 1_B))=\Sigma_j(\overline{\hat G^{-1/2}\xi_j}\otimes a^*_j).
$$
Then the image of the element $(X\otimes\overline\xi)^*=p(X^\bullet\otimes\overline{\xi^\bullet})=p(X^\bullet\otimes\overline{\overline{\hat
G^{1/2}\xi}})\in A_F$ equals to
$$
\Sigma_j(\overline{\hat G^{1/2}\xi},\overline{\hat G^{-1/2}\xi_j})a^*_j=(\Sigma_j(\xi_j,\xi)a_j)^*,
$$
which shows that the homomorphism is $*$-preserving. Passing to the $C^*$-completion, we have the first part of the proof.

Conversely, let us start with a unitary weak tensor functor $F$, construct a unital $\mathfrak G$-$C^*$-algebra $(A_F,{\mathfrak a}_F)$, and consider the spectral functor $F'$ associated with it. For any irreducible $U^x\in UCorep(\mathfrak G),\ x\in\Omega$,
fix an orthonrmal basis $\{\xi_i\}\in H_x$, then the space $F'(U_x)$ consists of vectors of the form $\Sigma_i(\xi_i\otimes X\otimes\overline \xi_i)$, where $X=F(U^x)$. The map $X\mapsto\Sigma_i(\xi_i\otimes X\otimes\overline \xi_i)$ from $F(U^x)$ to $F'(U^x)$ is clearly $A_\varepsilon $-bilinear, let us check that it is isometric. Taking $X'=\Sigma_i(\xi_i\otimes X\otimes\overline \xi_i),\ Y'=\Sigma_i(\xi_i\otimes Y\otimes\overline \xi_i)$ in $F'(U^x)$, we compute:
$$
<X',Y'>=\Sigma_i(X\otimes\overline \xi_i)^*(Y\otimes\overline \xi_i)=
$$
$$
=p(\Sigma_i(X^\bullet\otimes\overline{\overline{\hat G^{1/2}\xi_i}})(Y\otimes\overline \xi_i)=p(J(X^\bullet\otimes Y)\otimes\overline{R_{U^x}(1_B)}).
$$
Lemma \ref{bullet} and the fact that the morphism $R_{U^x}:B_s\to\overline{U^x}\otop U^x$ is an isometry imply that the last expression equals to $<X,Y>$, so the isomorphisms $F(U^x)\cong F'(U^x)$ are unitary and extend uniquely to a natural unitary isomorphism of functors $F$ and $F'$. Finally, one can check directly that this isomorphism is monoidal.
\hfill\end{dm}

\begin{proposition} \label{eqmod}
Let $\mathfrak G$ be a regular coconnected finite quantum groupoid and $\mathcal M$ be a strict right $UCorep(\mathfrak G)$-module $C^*$-category with 
generator $M$. If $(A,\mathfrak a)$ is a unital $\mathfrak G$-$C^*$-algebra constructed by this data in Lemma \ref{co}, then the category $\mathcal D_A$ 
(see Definition \ref{pluton}) is unitarily equivalent, as a $UCorep(\mathfrak G)$-module $C^*$-category, to $\mathcal M$, via an equivalence sending $A$ 
to $M$.
\end{proposition}

\begin{dm} As we have seen, $(F,J)$ is a weak tensor functor. Note that there are canonical isomorphisms of vector spaces
$$
F(U)\cong \mathcal D_A(A,U\otimes_{B_s} A)
$$
that map $\Sigma_i(\xi_i\otimes_{B_s} a_i)\in F(U)$ into the morphism $a\mapsto\Sigma_i(\xi_i\otimes_{B_s} a_i a)$. Therefore,
the spectral functor is naturally unitarily monoidally isomorphic to the weak tensor functor $F':UCorep(\mathfrak G)\to Corr(R)$ defined by $\mathcal D_A$ 
as in Proposition \ref{modcatweak}, where $R=End(A)$. If $\psi:F'\to F$ is such an isomorphism, then $\psi:A=F'(U^\varepsilon )\to F(U^\varepsilon )=A$ 
is the identity map since it is a bimodule map such that $\psi\circ J=J'(\psi\otimes\psi)$.

Let us now define a functor of linear categories $E:\tilde D_A\to\tilde{\mathcal M}$, where $\tilde D_A\subset D_A$ and $\tilde{\mathcal M}\subset\mathcal M$ 
are full subcategories consisting of objects $U\otimes_{B_s} A$ and $U\boxtimes M$,
respectively. We put $E(U\otimes_{B_s} A)=U\boxtimes M$ on objects and $E(T)=\psi(T)$ on morphisms $T\in D_A(A,U\otimes_{B_s} A)$.
More generally, if $T\in D_A(U\otimes_{B_s} A,V\otimes_{B_s} A)$, where $U,V\in UCorep(\mathfrak G)$, then $(id_{\overline U}
\otimes T)(R_U\otimes id_A)\in D_A(A,\overline U\otop V\otimes_{B_s} A)$ is Frobenius reciprocity isomorphism with inverse sending $S\in D_A(A,
\overline U\otop V\otimes_{B_s} A)$ to $(\overline R^*_U\otimes id\otimes id)(id_U\otimes S)$. We can define
similar isomorphisms in $\mathcal M$ and then define linear isomorphisms
$$
E:\mathcal D_A(U\otimes_{B_s} A,V\otimes_{B_s} A)\to\mathcal M(U\otimes_{B_s} M,V\otimes_{B_s} M)
$$
by $E(T)=(\overline R^*_U\otimes id\otimes id)[id_U\otimes\psi((id_{\overline U}\otimes T)(R_U\otimes id_A))]$.

Let us note that the naturality of $\psi$ implies that if $T:U\otimes_{B_s} A\to V\otimes_{B_s} A$, $S:V\to W$, where $U,V,W\in UCorep(\mathfrak G)$,
then:
\begin{equation} \label{natur}
E(id_W\otimes T)=id_W\otimes E(T)\quad\text{and}\quad E((S\otimes id)T)=(S\otimes id)E(T).
\end{equation}
Consider now morphisms $Q:U\otimes_{B_s} A\to V\otimes_{B_s} A$ and $T:V\otimes_{B_s} A\to W\otimes_{B_s} A$, and define the morphisms
$P=(id_{\overline U}\otimes Q)(R_U\otimes id_A):A\to(\overline U\otop V)\otimes_{B_s} A$ and $S=(id_{\overline V}\otimes T)(R_V\otimes id_A):
A\to(\overline V\otop W)\otimes_{B_s} A$, which give:
$$
TQ=(\overline R^*_V\otimes id_W\otimes id_A)(id_V\otimes S)(\overline R^*_U\otimes id_V\otimes id_A)(id_U\otimes P)=
$$
$$
=(\overline R^*_V\otimes id_W\otimes id_A)(\overline R^*_U\otimes id_V\otimes id_{\overline V}\otimes id_W\otimes id_A)
(id_U\otimes id_{\overline U}\otimes S)(id_U\otimes P)=
$$
$$
=(\overline R^*_U\otimes \overline R^*_V\otimes id_W\otimes id_A)(id_U\otimes J'(P\otimes S)),
$$
where $J'(P\otimes S)=(id_{\overline U}\otimes id_V\otimes S)P:A\to\overline U\otimes V\otimes\overline V\otimes W\otimes A$.
A similar calculation gives
$$
E(T)E(Q)=(\overline R^*_U\otimes \overline R^*_V\otimes id_W\otimes id_M)(id_U\otimes J(\psi(P)\otimes\psi(S))),
$$
from where, using (\ref{natur}) and monoidality of $\psi$, we get $E(TQ)=E(T)E(Q)$, which means that $\psi J'(P\otimes S)=
J(\psi(P)\otimes\psi(S))$. Therefore, $E$ is a functor, and since it is surjective on objects and fully faithful, it is an
equivalence of linear categories $\tilde{\mathcal D}_A$ and $\tilde M$.

Next, let us show that $E$ is unitary, i.e., $E(T^*)=E(T)^*$ on morphisms. First, let $T:A\to U\otimes_{B_s} A$. Since $\psi$
is unitary and $\psi|_A=id_A$, we have for any $S:A\to U\otimes_{B_s} A$:
$$
E(T)^*E(S)=\psi(T)^*\psi(S)=<\psi(T),\psi(S)>=<T,S>=
$$
$$
=T^*S=E(T^*S)=E(T^*)E(S).
$$
As $S$ is arbitrary, this implies that $E(T^*)=E(T)^*$, and using (\ref{natur}), we also have $E((T\otimes id)^*)=E(T\otimes id)^*$.
But any morphism in $\tilde{\mathcal D}_A$ is a composition of two morphisms: one of the above form $T\otimes id_V$ and another of
the form $id_M\otimes S$ for some morphism $S$ in $UCorep(\mathfrak G)$. As a consequence of (\ref{natur}), we have
$E(id_M\otimes S)^*=(id_M\otimes S)^*=E((id_M\otimes S)^*)$, it follows that $E$ is unitary.

Further, if we define $J=J_{U\otimes A,V}:V\otimes_{B_s} E(U\otimes_{B_s} A)\to E((V\otop U)\otimes_{B_s} A)$ to be the identity
maps, the relations (\ref{natur}) show that we get a natural isomorphism of bilinear functors $\cdot\otimes E(\cdot)$ and
$E(\cdot\otimes\cdot)$. Therefore, $(E,J)$ is a unitary equivalence of $UCorep(\mathfrak G)$-$C^*$-module categories $\tilde{\mathcal D}_A$
and $\tilde{\mathcal M}$.

Finally, since $\mathcal D_A$ and $\mathcal M$ are completions of these categories with respect to subobjects, the equivalence between
$\tilde{\mathcal D}_A$ and $\tilde{\mathcal M}$ extends uniquely, up to a natural unitary isomorphism, to a unitary equivalence between
the $UCorep(\mathfrak G)$-$C^*$-module categories $\mathcal D_A$ and $\mathcal M$.
\hfill\end{dm}

Now we are ready to prove Theorem \ref{main}.

\begin{dm} Due to the previous proposition, it remains to show that two unital $\mathfrak G$-$C^*$-algebras, $(A_1,\mathfrak a_1)$ and
$(A_2,\mathfrak a_2)$, are isomorphic if and only if the pairs $(\mathcal D_{A_1},A_1)$ and $(\mathcal D_{A_2},A_2)$ are unitarily equivalent.

First, given such equivalent pairs, we have the isomorphism of the corresponding spectral subalgebras $(A_1)_\varepsilon =End(M_1)$ and
$(A_2)_\varepsilon =End(M_2)$. Identifying the above algebras via this isomorphism, we have a natural unitary monoidal isomorphism of the weak tensor functors constructed in Proposition \ref{modcatweak} which implies a natural unitary monoidal isomorphism of the corresponding spectral functors. Now theorem \ref{spec-weak} gives the needed isomorphism of unital $\mathfrak G$-$C^*$-algebras. Conversely, isomorphic unital $\mathfrak G$-$C^*$-algebras clearly produce unitarily equivalent classes of pairs of the form
$(\mathcal M, M)$.
\hfill\end{dm}

Note that: (i) one can precise the definition of the equivalence of module functors between pairs $(\mathcal M,M)$ as in  \cite{CY},
Theorem 6.4; (ii) under the above equivalence, the unital $C^*$-algebra $A_\varepsilon $ is isomorphic to $End_{\mathcal M}(M)\otimes B_s$.

\begin{corollary} \label{abstract}
Let $\mathcal M$ be a strict left module $C^*$-category over a strict rigid finite $C^*$-tensor category $\mathcal C$, $M$ be a generator in $\mathcal M$, 
and denote by $R$ the unital $C^*$-algebra $End(M)$. Then there exist a regular biconnected finite quantum groupoid $\mathfrak G$ (even with commutative 
base) and a unital $\mathfrak G$-$C^*$-algebra $(A,\mathfrak a)$ such that $\mathcal C$ is equivalent to $UCorep(\mathfrak G)$ as $C^*$-tensor categories 
and $\mathcal M$ is equivalent to $\mathcal D_A$ as left $UCorep(\mathfrak G)$-module $C^*$-categories via an equivalence that maps $M$ to $A$.
\end{corollary}

Indeed, the existence of $\mathfrak G$ is guaranteed by Theorem \ref{reconstruction}, and the second statement - by Proposition \ref{eqmod}.

\begin{corollary} \label{simperg} If $\mathfrak G$ is regular and coconnected, then $A_\varepsilon =A^{\mathfrak a}\alpha(B_s)$.

Indeed, we have seen that $A_\varepsilon = End_{\mathcal M}(M)\otimes B_s$ and that $A^{\mathfrak a}= End_{\mathcal M}(M)$.
\end{corollary}

\begin{example} \label{regular} 
The  $C^*$-algebra $B$ with coproduct $\Delta$ viewed as $\mathfrak G$-$C^*$-algebra, corresponds to the $UCorep(\mathfrak G)$-module $C^*$-category
$\mathcal Corr_f(B_s)$ with generator $M=B_s$: for any element $U\in UCorep(\mathfrak G)$ and $N\in Corr_f(B_s)$, one defines $U\boxtimes N:=
F(U)\otimes_{B_s} N$, where the functor $F: UCorep(\mathfrak G)\to Corr_f(B_s)\ (F(U)=H_U)$ is the forgetful functor. Indeed, if one identifies
$\mathcal M(B_s,H_U)$ with $H_U$, we get an isomorphism of the algebra $\tilde A$ constructed from the pair $(\mathcal M,M)$ onto
$\tilde B=\underset{U}\bigoplus(H_U\otimes\overline H_U)$ and then an isomorphism $A\cong B=\underset{x\in\hat G}\bigoplus (H_x\otimes\overline H_x)$ 
such that $p:\tilde A\to A$ turns into the map $\tilde B\to B$ sending $\xi\otimes\overline\eta\in H_U\otimes\overline H_U$ into the matrix coefficient $U_{\xi,\eta}$.
\end{example}
\end{section}
%%%%%%%%%%%%%%%%%%%%%%%%%%%%%%%%%%%%%%%%%%%%%%%%%%%%%%%%%%%%%%%%%%%%%%%%%%%%%%%%%%%%%%%%%%%%%%%%%%%%%%%%%%%%%%%%%%%%%%%

\bibliographystyle{plain}
\bibliography{biblio1}

\begin{thebibliography}{10}

\bibitem{BS1}
Saad Baaj and Georges Skandalis.
\newblock {$C^\ast$}-alg\`ebres de {H}opf et th\'eorie de {K}asparov
  \'equivariante.
\newblock {\em $K$-Theory}, 2(6):683--721, 1989.

\bibitem{Boca}
Florin~P. Boca.
\newblock Ergodic actions of compact matrix pseudogroups on {$C^*$}-algebras.
\newblock {\em Ast\'erisque}, (232):93--109, 1995.
\newblock Recent advances in operator algebras (Orl{\'e}ans, 1992).

\bibitem{BNSz}
Gabriella B{\"o}hm, Florian Nill, and Korn{\'e}l Szlach{\'a}nyi.
\newblock Weak {H}opf algebras. {I}. {I}ntegral theory and {$C^*$}-structure.
\newblock {\em J. Algebra}, 221(2):385--438, 1999.

\bibitem{BoSz}
Gabriella B{\"o}hm and Korn{\'e}l Szlach{\'a}nyi.
\newblock Weak {$C^\ast$}-{H}opf algebras: the coassociative symmetry of
  non-integral dimensions.
\newblock In {\em Quantum groups and quantum spaces ({W}arsaw, 1995)},
  volume~40 of {\em Banach Center Publ.}, pages 9--19. Polish Acad. Sci.,
  Warsaw, 1997.

\bibitem{BSz}
Gabriella B{\"o}hm and Korn{\'e}l Szlach{\'a}nyi.
\newblock Weak {H}opf algebras. {II}. {R}epresentation theory, dimensions, and
  the {M}arkov trace.
\newblock {\em J. Algebra}, 233(1):156--212, 2000.

\bibitem{CY}
Kenny De~Commer and Makoto Yamashita.
\newblock Tannaka-{K}re\u\i n duality for compact quantum homogeneous spaces.
  {I}. {G}eneral theory.
\newblock {\em Theory Appl. Categ.}, 28:No. 31, 1099--1138, 2013.

\bibitem{CY1}
Kenny De~Commer and Makoto Yamashita.
\newblock Tannaka-{K}re\u\i n duality for compact quantum homogeneous spaces
  {II}. {C}lassification of quantum homogeneous spaces for quantum {$\rm
  SU(2)$}.
\newblock {\em J. Reine Angew. Math.}, 708:143--171, 2015.

\bibitem{E2}
Michel Enock.
\newblock Measured quantum groupoids in action.
\newblock {\em M\'em. Soc. Math. Fr. (N.S.)}, (114):ii+150 pp. (2009), 2008.

\bibitem{EGNO}
Pavel Etingof, Shlomo Gelaki, Dmitri Nikshych, and Victor Ostrik.
\newblock {\em Tensor categories}, volume 205 of {\em Mathematical Surveys and
  Monographs}.
\newblock American Mathematical Society, Providence, RI, 2015.

\bibitem{Ha}
T~Hayashi.
\newblock A canonical tannaka duality for semi finite tensor categories.
\newblock Preprint, math.QA/9904073, 1999.

\bibitem{L}
E.~C. Lance.
\newblock {\em Hilbert {$C^*$}-modules}, volume 210 of {\em London Mathematical
  Society Lecture Note Series}.
\newblock Cambridge University Press, Cambridge, 1995.
\newblock A toolkit for operator algebraists.

\bibitem{ML}
Saunders Mac~Lane.
\newblock {\em Categories for the working mathematician}, volume~5 of {\em
  Graduate Texts in Mathematics}.
\newblock Springer-Verlag, New York, second edition, 1998.

\bibitem{M}
Camille Mevel.
\newblock Exemples et applications des groupoides quantiques finis,
  https://tel.archives-ouvertes.fr/tel-00498884/document.
\newblock Th\`ese, Universit\'e de Caen, 2010.

\bibitem{NY}
S~Neshveyev and M~Yamashita.
\newblock Categorical duality for {Y}etter-{D}rinfeld algebras.
\newblock Preprint, arXiv: 1310.4407v4 [math.OA], 2013.

\bibitem{Nes1}
Sergey Neshveyev.
\newblock Duality theory for nonergodic actions.
\newblock {\em M\"unster J. Math.}, 7(2):413--437, 2014.

\bibitem{NT}
Sergey Neshveyev and Lars Tuset.
\newblock Hopf algebra equivariant cyclic cohomology, {$K$}-theory and index
  formulas.
\newblock {\em $K$-Theory}, 31(4):357--378, 2004.

\bibitem{Nes}
Sergey Neshveyev and Lars Tuset.
\newblock {\em Compact quantum groups and their representation categories},
  volume~20 of {\em Cours Sp\'ecialis\'es [Specialized Courses]}.
\newblock Soci\'et\'e Math\'ematique de France, Paris, 2013.

\bibitem{NY1}
Sergey Neshveyev and Makoto Yamashita.
\newblock Categorical duality for {Y}etter-{D}rinfeld algebras.
\newblock {\em Doc. Math.}, 19:1105--1139, 2014.

\bibitem{NTV}
Dmitri Nikshych, Vladimir Turaev, and Leonid Vainerman.
\newblock Invariants of knots and 3-manifolds from quantum groupoids.
\newblock In {\em Proceedings of the {P}acific {I}nstitute for the
  {M}athematical {S}ciences {W}orkshop ``{I}nvariants of {T}hree-{M}anifolds''
  ({C}algary, {AB}, 1999)}, volume 127 1-2, pages 91--123, 2003.

\bibitem{NV5}
Dmitri Nikshych and Leonid Vainerman.
\newblock Algebraic versions of a finite-dimensional quantum groupoid.
\newblock In {\em Hopf algebras and quantum groups ({B}russels, 1998)}, volume
  209 of {\em Lecture Notes in Pure and Appl. Math.}, pages 189--220. Dekker,
  New York, 2000.

\bibitem{NV1}
Dmitri Nikshych and Leonid Vainerman.
\newblock A characterization of depth 2 subfactors of {${\rm II}_1$} factors.
\newblock {\em J. Funct. Anal.}, 171(2):278--307, 2000.

\bibitem{NV2}
Dmitri Nikshych and Leonid Vainerman.
\newblock A {G}alois correspondence for {${\rm II}_1$} factors and quantum
  groupoids.
\newblock {\em J. Funct. Anal.}, 178(1):113--142, 2000.

\bibitem{NV}
Dmitri Nikshych and Leonid Vainerman.
\newblock Finite quantum groupoids and their applications.
\newblock In {\em New directions in {H}opf algebras}, volume~43 of {\em Math.
  Sci. Res. Inst. Publ.}, pages 211--262. Cambridge Univ. Press, Cambridge,
  2002.

\bibitem{Ni}
Florian Nill.
\newblock Axioms for weak bialgebras.
\newblock Preprint,arXiv: math/9805104 [math.QA], 1998.

\bibitem{Ped}
Gert~K. Pedersen.
\newblock {\em {$C^{\ast} $}-algebras and their automorphism groups}, volume~14
  of {\em London Mathematical Society Monographs}.
\newblock Academic Press, Inc. [Harcourt Brace Jovanovich, Publishers],
  London-New York, 1979.

\bibitem{Pf1}
Hendryk Pfeiffer.
\newblock Finitely semisimple spherical categories and modular categories are
  self-dual.
\newblock {\em Adv. Math.}, 221(5):1608--1652, 2009.

\bibitem{SZ}
K.~Szlach\'anyi.
\newblock Finite quantum groupoids and inclusions of finite type.
\newblock {\em Fields Inst.Commun.}, 30:314--343. AMS Providence, RI, 2001.

\bibitem{TY}
Daisuke Tambara and Shigeru Yamagami.
\newblock Tensor categories with fusion rules of self-duality for finite
  abelian groups.
\newblock {\em J. Algebra}, 209(2):692--707, 1998.

\bibitem{VV}
Leonid Vainerman and Jean-Michel Vallin.
\newblock On {\Large z}/2{\Large z}-extensions of pointed fusion categories.
\newblock {\em Banach Center Publications}, 98:343--366, 2012.

\bibitem{Val1}
Jean-Michel Vallin.
\newblock Groupo\"\i des quantiques finis.
\newblock {\em J. Algebra}, 239(1):215--261, 2001.

\bibitem{Val5}
Jean-Michel Vallin.
\newblock Multiplicative partial isometries and finite quantum groupoids.
\newblock In {\em Locally compact quantum groups and groupoids ({S}trasbourg,
  2002)}, volume~2 of {\em IRMA Lect. Math. Theor. Phys.}, pages 189--227. de
  Gruyter, Berlin, 2003.

\bibitem{Val2}
Jean-Michel Vallin.
\newblock Actions and coactions of finite quantum groupoids on von {N}eumann
  algebras, extensions of the matched pair procedure.
\newblock {\em J. Algebra}, 314(2):789--816, 2007.

\bibitem{Wor}
S.~L. Woronowicz.
\newblock Tannaka-{K}re\u\i n duality for compact matrix pseudogroups.
  {T}wisted {${\rm SU}(N)$} groups.
\newblock {\em Invent. Math.}, 93(1):35--76, 1988.

\end{thebibliography}

\end{document}